\documentclass[preprint,12pt]{elsarticle}
\usepackage[top=1.0in,left=0.80in,bottom=1.0in,right=0.80in]{geometry}

\makeatletter
\def\ps@pprintTitle{%
 \let\@oddhead\@empty
 \let\@evenhead\@empty
 \def\@oddfoot{}%
 \let\@evenfoot\@oddfoot}
\makeatother
\usepackage{url}
\usepackage{amssymb}
\usepackage{latexsym}
\usepackage{mathtools}
\usepackage{mathrsfs}
\usepackage{amsfonts}
\usepackage{multirow}
\usepackage{color}
\usepackage{colortbl}
\usepackage{amsthm}

\usepackage[ruled,vlined]{algorithm2e}

\newtheorem{theorem}{Theorem}
\newtheorem{remark}{Remark}

\newtheorem{corollary}{Corollary}

\newtheorem{lemma}{Lemma}
\newtheorem{assumption}{Assumption}

\DeclareMathOperator*{\argmin}{arg\,min}

\newcommand{\cG}{\mathcal G}
\newcommand{\cH}{\mathcal H}
\newcommand{\co}{\textrm {Cov}}
\newcommand{\ou}{\overline u}
\newcommand{\C}{\mathcal C}

\newcommand{\E}{\mathbb E}
\newcommand{\n}{\mathbf{n}
\newcommand{\Ga}{\Upsilon}}
\newcommand\norm[1]{\left\lVert#1\right\rVert}

\journal{}

\begin{document}

\begin{frontmatter}



\title{{\LARGE \bf Adaptive regularisation for ensemble Kalman inversion}}

 \author[label1]{{\bf Marco Iglesias }}
 \ead{marco.iglesias@nottingham.ac.uk}
  \author[label1]{{\bf Yuchen Yang }}

\address[label1]{School of Mathematical Sciences, University of Nottingham, Nottingham, NG72RD}


\address{}

\begin{abstract}
We propose a new regularisation strategy for the classical ensemble Kalman inversion (EKI) framework. The strategy consists of: (i) an adaptive choice for the regularisation parameter in the update formula in EKI, and (ii) criteria for the early stopping of the scheme. In contrast to existing approaches, our parameter choice does not rely on additional tuning parameters which often have severe effects on the efficiency of EKI. We motivate our approach using the interpretation of EKI as a Gaussian approximation in the Bayesian tempering setting for inverse problems.  We show that our parameter choice controls the symmetrised Kulback-Leibler divergence between consecutive tempering measures. We further motivate our choice using a heuristic statistical discrepancy principle.

We test our framework using electrical impedance tomography with the complete electrode model. Parameterisations of the unknown conductivity are employed which enable us to characterise both smooth or a discontinuous (piecewise-constant) fields. We show numerically that the proposed regularisation of EKI can produce efficient, robust and accurate estimates, even for the discontinuous case which tends to require larger ensembles and more iterations to converge. We compare the proposed technique with a standard method of choice and demonstrate that the proposed method is a viable choice to address computational efficiency of EKI in practical/operational settings.


\end{abstract}







\end{frontmatter}


\section{Introduction}

Ensemble Kalman Inversion (EKI)  \cite{ILS13a,EnsembleYo,Neil,Tong1, Tong2} is a collection of algorithms that import ideas from the ensemble Kalman Filter (EnKF) developed in \cite{doi:10.1029/94JC00572}. While EnKF was devised for assimilating data into models for numerical weather prediction and oceanography, EKI has been applied to solve parameter identification problems arising from multiple disciplines. The initial versions of EKI were proposed for the calibration of oil-reservoir models \cite{EME2,Oliv}, and then transferred in \cite{ILS13a, EnsembleYo} to more generic PDE-constrained inverse problems settings. The development of EKI as an iterative solver for parameter identification problems has lead to numerous application including the calibration of climate  \cite{doi:10.1002/2017GL076101}, turbulent flow \cite{XIAO2016115}, finance \cite{1930-8337_2017_5_799} and biomedical \cite{Chido} models. EKI has been used for imaging and non-destructive testing including electrical impedance tomography \cite{Calvetti_2020}, seismic inversion \cite{seismic}, characterisation of thermophysical properties of walls \cite{DESIMON2018220} and composite materials \cite{Iglesias_2018}. More recently, EKI has also been applied for the solution of machine learning tasks \cite{ML}.


There is clear promise for the potential use of EKI as a practical and operational tool for PDE-constrained calibration and imaging arising from multiple applications in science and engineering. However, most EKI algorithms still rely on the appropriate selection of user-defined parameters that control the stability, accuracy and computational efficiency of EKI. Unfortunately, the lack of a general theory for the convergence of EKI means that there is no principled approach to select these parameter in an optimal fashion. The seamless version of EKI proposed in \cite{SS17} and further developed in \cite{SS18,Bl_mker_2019,Alfredo1,doi:10.1137/17M1132367} has provided enormous theoretical advances for understanding EKI within the more general frameworks of Stochastic Differential Equations (SDEs). However, to the best of our knowledge, the selection of those crucial EKI parameters are still chosen empirically. Among those parameters, the choice of the regularisation (inflation) parameter in EKI, or alternatively,  the step size from the discretisation of the SDE formulation of EKI often depends on additional tuning parameters which can only be informed after problem-specific numerical testing which maybe computationally intensive. 

The aim of this work is to introduce a simple, yet computationally efficient, regularisation strategy for EKI that does not rely on further tuning parameters. Our main focus is large-scale/high resolution imaging/identification settings in which there are two fundamental challenges: (i) the forward model is computationally costly and (ii) the unknown must be parameterised in a highly-complex manner so that key features from the truth can be extracted via EKI. Both challenges are intertwined since the latter means that EKI requires large ensembles and more iterations to achieve accurate identifications; thus the total computational cost of EKI can become unfeasible. An efficient and robust regularisation strategy within EKI is a key requirement in these practical settings. 

\subsection{The inverse problem framework with ensemble Kalman inversion.}
We work on a generic setting where the properties that we wish to identify are functions which, in turn, play the role of inputs (e.g. coefficients) of Partial Different Equations (PDEs) describing the underlying experiment/process. The observable variables in these experiment/processes are functionals of the solution of these PDEs. Under this inverse problem setting, the forward problem can be written in terms of a nonlinear operator $\mathcal{F}: \mathcal{X}\to \mathbb{R}^{M}$ that maps physical properties from an admissible space $\mathcal{X}$ to the space of observable quantities $\mathbb{R}^{M}$. Our work is focused on using EKI for solving the \textit{classical (deterministic) inverse problem} that consists of estimating the underlying physical property of the medium, denoted by $\kappa^{\dagger}$, given noisy measurements of $\mathcal{F}(\kappa^{\dagger})$, which we assume are given by 
\begin{eqnarray}\label{equ0}
y=\mathcal{F}(\kappa^{\dagger})+\eta
\end{eqnarray}
where $\eta$ is the unknown measurement error. We further assume $\eta$ is drawn from a centred Gaussian measure with known covariance $\Gamma$. In order to address this inverse problem, we define a class of suitable parameterisations  $\mathcal{P}:\mathcal{H}\to \mathcal{X}$ that enable us to characterise our estimate of the unknown as $\kappa=\mathcal{P}(u)$, i.e. in terms of an input (function) $u\in \mathcal{H}$ which we calibrate within the EKI framework. We pose the parameterised (in terms of $u$) inversion via 
\begin{eqnarray}\label{equ1}
u^{\ast}=\argmin_{u\in \mathcal{S}_{0}}\Phi (u;y),
\end{eqnarray}
where $\mathcal{S}_{0}\subset \mathcal{H}$ is a user-defined subspace of admissible solutions and $\Phi (u;y)$ is the functional defined by
\begin{eqnarray}\label{equ1B}
\Phi (u;y)\equiv \frac{1}{2}\norm{\Gamma^{-1/2}( y-\cG(u) )}^2,
\end{eqnarray}
where $\cG=\mathcal{F}\circ\mathcal{P}$ is the forward (parameter-to-output) map and $\norm{\cdot}$ denotes the $M$-dimensional Euclidean norm. We use EKI to approximate (\ref{equ1}) and to provide an estimate of $\kappa^{\dagger}$ via $\kappa^{\ast}=\mathcal{P}(u^{\ast})$. 


For the PDE-constrained identification problems that we wish to address, $\mathcal{\cG}$ is often a compact operator which leads to the ill-posedness of (\ref{equ1}) in the sense of stability \cite{Engl,Kaipio:1338003}. Although numerous regularisation techniques \cite{Engl,Iterative} can be used to address ill-posedness, most of them require the computation of the Frech\'{e}t derivative of the forward map, as well as the corresponding adjoint operator. This constitutes a substantial limitation in many practical applications where the map $\mathcal{F}$  is only accessible via commercial software simulations with no adjoint built-in functionalities. Such a practical limitation has given rise to a large body of work on EKI techniques that stem from EnKF \cite{doi:10.1029/94JC00572} and which, in turn, do not require derivatives of the forward map.

While there are several versions of EKI algorithms to approximate (\ref{equ1}), here we focus on the classical, perturbed-observation EKI displayed in Algorithm \ref{Al1}. This generic version of EKI involves selecting an initial ensemble of $J$ particles $\{u_{0}^{(j)}\}_{j=1}^{J}\subset \mathcal{H}$. Then, each particle is iteratively updated according to the following formula
\begin{eqnarray}\label{equ2}
u_{n+1}^{(j)}=u_{n}^{(j)}+\C_{n}^{u\cG}(\C_{n}^{\cG \cG} +\alpha_{n}\Gamma  )^{-1}(y+\sqrt{\alpha_n}\xi_{n}-\cG(u_{n}^{(j)}))
\end{eqnarray}
where $\alpha_{n}$ is a tuning (regularisation) parameter, $\xi_{n}\sim N(0,\Gamma)$ is perturbation of the data, and $\C_{n}^{u\cG}$ and $\C_{n}^{\cG \cG}$ are empirical covariances defined in (\ref{equ3})-(\ref{equ3B}). The running estimate of the unknown is obtained via the ensemble mean
\begin{eqnarray}\label{equ4}
 \overline{u}_{n+1}\equiv \frac{1}{J}\sum_{j=1}^{J} u_{n+1}^{(j)}.
\end{eqnarray}

Using informal arguments, it can be shown (see \ref{Ape_LM} and the work in \cite{EnsembleYo}) that the ensemble $ \overline{u}_{n+1}$ can be seen as an ensemble approximation of 
\begin{eqnarray}\label{eqA:14}
m_{n+1}= \argmin_{u\in \mathcal{S}_{0}}\Big\{\frac{1}{2}\norm{ \Gamma^{-1/2}(y-\cG(m_n)-D\cG(m_n)(u-m_{n})}^2+\frac{\alpha_{n}}{2}\norm{ \C_{n}^{-1/2} (u-m_{n})}_{\cH}^{2}\Big\}
\end{eqnarray}
where $D\cG$ denotes the Frech\'{e}t  derivative of $\cG$, $\C_{n}$ is a covariance operator that we define in \ref{Ape_LM}, and $\mathcal{S}_{0}\equiv \text{span}\{u_{0}^{(j)}\}_{j=1}^{J}$ is the subspace generated by the initial ensemble. If $\C_{n}$ is the identity operator, (\ref{eqA:14}) is the standard Levenberg-Marquardt (LM) scheme \cite{Mar} applied for the computation of (\ref{equ1}). Note that (\ref{eqA:14}) can also be interpreted as iterative Tikhonov regularisation applied to the linearisation of $\mathcal{G}$.

The link between EKI and (\ref{eqA:14}) is very useful because (i) it motivates EKI as a derivative-free solver for (\ref{equ1}) and (ii) it reveals the role of $\alpha_{n}$ in (\ref{equ2}) as a Tikhonov regularisation parameter. According to the theory in \cite{Hanke}, $\alpha_{n}$ must be carefully selected, together with the stopping criteria, in order to ensure the stability of the LM scheme. The approach for selecting $\alpha_{n}$ in the LM proposed in \cite{Hanke} has been adapted to the EKI framework in \cite{EnsembleYo,Iglesias_2018}, and subsequently used in \cite{pub:28949,Neil,pub:29192,seismic,Chido,cited1}. As we discuss in the next section, this approach relies on tuning parameters that, unless carefully chosen, can lead to unnecessary large number of iterations $n^*$. Since the main computational cost of Algorithm \ref{Al1} is $n^*{}\times J$,  it is clear that a large $n^*$ is detrimental to the computational efficiency of EKI.

\subsection{EKI as a Gaussian approximation in the Bayesian tempering setting}

Although the goal of most of the existing applications of EKI is to solve the deterministic problem in (\ref{equ1}), the role of EKI within the Bayesian setting for parameter identification can be useful for identifying suitable choices of the regularisation parameter $\alpha_{n}$. In the Bayesian setting we put a prior measure, $\mu_{0}(u)=\mathbb{P}(u)$, on the unknown $u$ that we wish to infer. Given measurements, $y$, the Bayesian inverse problem consists of approximating the posterior $\mu(u)\equiv \mathbb{P}(u\vert y)$ which, from Bayes' rule \cite{Andrew} is given by 
\begin{eqnarray}\label{eqA:4}
\mu(du) \propto   \mu_{0}(du)\exp\Big[-\frac{1}{2}\norm{ \Gamma^{-1/2}(y-\cG(u)) }^2\Big],
\end{eqnarray}
where we have made the standard assumption that $y=\cG(u)+\eta$ with $\eta\sim N(0,\Gamma)$. Modern computational approaches \cite{BJMS15,Kantas,LATZ2018154,Ruchi_2019,tawn_roberts_2019} for high-dimensional Bayesian inverse problems, use the tempering approach that consists of introducing $N$ intermediate measures $\{\mu_{t_{n}}\}_{n=1}^{N}$ between the prior and the posterior. These measures are defined by
\begin{eqnarray}\label{eqA:5}
\mu_{t_n}(du)\propto  \mu_{t_{0}}(du)\exp\Big[-\frac{t_{n}}{2}\norm{ \Gamma^{-1/2}(y-\cG(u)) }^2\Big],
\end{eqnarray}
 where $\{t_{n}\}_{n=1}^{N}$ are tempering parameters that satisfy:
\begin{eqnarray}\label{eqA:6B}
t_{0}\equiv 0<t_{1}<t_{2}<\dots<t_{N}<t_{N+1}\equiv 1.
\end{eqnarray}
Note that $n=0$ and $n=N+1$ in (\ref{eqA:5}) yields the prior ($\mu_{t_{0}}=\mu_{0}$) and posterior ($\mu_{t_{N+1}}=\mu$), respectively. From expression (\ref{eqA:5}) we obtain the following recursive formula for the intermediate measures:
\begin{eqnarray}\label{eqA:7}
\frac{ \mu_{t_{n+1}}(du)}{ \mu_{t_{n}}(du)}\propto \exp\Big[-\frac{1}{2}\norm{ (\alpha_{n}\Gamma)^{-1/2}(y-\cG(u))}^2\Big],
\end{eqnarray}
where 
\begin{eqnarray}\label{eqA:8}
\alpha_{n}^{-1}=t_{n+1}-t_{n}.
\end{eqnarray}
From (\ref{eqA:6B}) it follows that
\begin{eqnarray}\label{eqA:9}
\sum_{n=0}^{N}\alpha_{n}^{-1}=1.
\end{eqnarray}

In (\ref{eqA:8}) we employ the same notation that we use for the regularisation parameter in EKI (see eq. (\ref{equ2})), because the ensemble computed at the $n$th iteration of EKI, is an ensemble approximation of a Gaussian measure that, in turn, approximates\footnote{Even when $J$ is large, these approximations are only exact in the linear-Gaussian case (i.e. $\cG$ linear and $\mu_{0}$ Gaussian)} the intermediate distribution $\mu_{t_n}$ in (\ref{eqA:5}) (see \cite{Iglesias_2018} and \ref{derivative_free}). Therefore, in the Bayesian setting, controlling the regularisation parameter $\alpha_{n}$ in EKI means to gradually transition between prior and posterior in order to facilitate more accurate sampling of the intermediate measures. 



The link between EKI and the Bayesian tempering setting has been recently explored in \cite{Iglesias_2018}, where the selection of $\alpha_{n}$ is borrowed from the adaptive-tempering Sequential Monte Carlo (SMC) method of \cite{Kantas}. In this approach, $\alpha_{n}$ is selected based on user-defined threshold to the Effective Sample Size (ESS) which, in SMC methods, is used to determine the quality in the population of particles that approximate each intermediate measure $\mu_{t_n}$. However, this selection of $\alpha_{n}$ involves the computation of the likelihood between consecutive measures (\ref{eqA:7}). When the number of observations is large and/or a the measurement noise is small, this likelihood can take very small values unless large $\alpha_{n}$'s are chosen. Consequently, many iterations may be be needed to satisfy condition (\ref{eqA:9}). Furthermore, the approach of \cite{Iglesias_2018} requires the aforementioned user-defined threshold on the (ESS) which may substantially affect the efficiency of EKI.

\subsection{Our contribution.} 

In this work we propose a novel adaptive choice of $\alpha_{n}$ in EKI that does not require any tuning parameters. This new selection of $\alpha_{n}$ depends on the number of observations, $M$, as well as the the values of the least-squares functional (\ref{equ1B}) for the ensemble of particles, i.e.
\begin{equation} \label{DMC}
\Phi_{n}\equiv \big\{\Phi(u_{n}^{(j)};y)\big\}_{j=1}^{J}=\Bigg\{\frac{1}{2}\norm{\Gamma^{-1/2}( y-\cG(u_{n}^{(j)}) )}^2\Bigg\}_{j=1}^{J}
\end{equation}
More specifically, we select $\alpha_{n}^{-1}$ via
\begin{equation} \label{eq: data-misfit controller}
\alpha_{n}^{-1} = \min \left\lbrace \max \left\lbrace  \frac{M}{ 2 \overline{\Phi}_{n}} ,\, \sqrt{        \frac{M}{2\sigma_{\Phi_{n}}^2 }     } \right\rbrace ,\, 1 -t_{n}  \right\rbrace
\end{equation}
where $\overline{\Phi}_{n}$ and $\sigma_{\Phi_{n}}^2 $ are the empirical mean and variance of $\Phi_{n}$, and 
\begin{equation} \label{eq:CCC}
\arraycolsep=1.4pt\def\arraystretch{1.7}
t_{n}=\left\{\begin{array}{cc}
\sum_{j=0}^{n-1}\alpha_{j}^{-1} &\quad \text{if}\quad n\geq 1,\\
0 & \quad \text{if}\quad  n=0.\end{array}\right.
\end{equation}
We stop the EKI algorithm at the iteration, $n^*$, such that
$$\alpha_{n^*}^{-1}=1-t_{n^*}$$
which satisfies constraint \eqref{eqA:9} from the tempering/annealing setting using $N=n^*$ intermediate measures. Noting that the selection of $\alpha_{n}^{-1}$ in \eqref{eq: data-misfit controller} is determined by the data misfit of the particles, we refer to our regularisation strategy as the \textit{data misfit controller} (DMC).

We motivate our DMC from the Bayesian perspective of EKI in which $\alpha_{n}^{-1}$ is the difference between consecutive tempering/annealing parameters (see eq. \eqref{eqA:8}). We show that the selection of $\alpha_{n}^{-1}$ in \eqref{eq: data-misfit controller} can be obtained by controlling, via imposing a threshold, the symmetrised Kulback-Leibler divergence (or Jeffreys' divergence) between two consecutive tempering measures. We then use a statistical discrepancy principle to select this threshold.

Jeffreys' divergence has been employed to design efficient selection of the step size between tempering measures when using MCMC to approximate the target posterior. For example, in the context of the \textit{tempered transitions} method, the work of \cite{tunning} uses Jeffreys divergence to select step sizes that increase acceptance rates. Jeffreys' divergence has also been used for path sampling: a method commonly used to estimate normalisation constants (e.g. for Bayesian model selection). The work of \cite{ps}, for example, shows that by controlling the Jeffreys' divergence between tempering measures, the error in path sampling estimator can be minimised. Both the works in \cite{tunning} and \cite{ps}, as well as various references therein,  employ continuous tempering in the finite-dimensional setting in order to characterise Jeffreys' divergence. Here we extend some of those results to our infinite-dimensional framework for EKI. 

We encode our regularisation strategy for EKI in an algorithm ({\bf EKI-DMC}) that we test via Electrical Impedance tomography (EIT) with the Complete Electrode Model (CEM). We show that the choice of $\alpha_{n}$ via the DMC that we import from the Bayesian formulation, leads to a stable, robust and computationally efficient EKI algorithm for solving classical (deterministic) inverse problems posed via (\ref{equ1}). We demonstrate the computational and practical advantages of {\bf EKI-DMC} over existing approaches in which $\alpha_{n}$ is borrowed from the theory for the LM scheme. We conduct test with two parameterisations of the unknown conductivity that allows to consider both continuous and piece-wise constant conductivities. In particular, we investigate the performance {\bf EKI-DMC} to infer regions characterised via a level-set function which is, in turn, parameterised with an anisotropic Whittle-Matern (WM) field.

In Section \ref{Liter} we further review the literature on existing approaches for the selection of $\alpha_{n}$ for the classical EKI framework. The proposed DMC (\ref{eq: data-misfit controller}) is motivated in Section \ref{Moti}, where we use continuous tempering to derive approximations to Jeffreys' divergence between tempering measures (subsection \ref{subsec: information gain}), as well as a statistical discrepancy principle to (subsection \ref{subsec: mean-variance pair}) for selecting the threshold on Jeffreys' divergence. In Section \ref{Num_testing} we provide a numerical investigation of the performance of the proposed {\bf EKI-DMC} algorithm for EIT. In Section \ref{Conclu} we discuss final remarks and conclusions. For completeness in \ref{Ape_LM} we motivate the classical EKI from the Bayesian tempering scheme, and from which we also draw links between the LM algorithm and EKI. Some technical proofs are included in \ref{app: appendix}.

\begin{algorithm}[h!]

\SetAlgoLined

\SetKwInOut{Input}{Input}
\Input{\begin{itemize}
\item[1)] $\{u_{0}^{(j)}\}_{j=1}^{J}$: Initial ensemble of inputs.
\item[2)] Measurements $y$ and covariance of measurement errors $\Gamma$.
\end{itemize}
}
Set $\{u_{n}^{(j)}\}_{j=1}^{J}=\{u_{0}^{(j)}\}_{j=1}^{J}$ and $\theta=0$\\
 \While{$\theta<1$ }{
\begin{itemize}
\item[(1)] Compute $\cG_{n}^{(j)}= \cG(u_{n}^{(j)}),\qquad j\in\{1,\dots, J\}$
\item[(2)] Compute regularisation parameter $\alpha_n$ and check for convergence criteria\newline

\vspace{-5mm}
 \textbf{if} converged \textbf{then} \newline
\vspace{-2.5mm}
\hspace{10mm} set $\theta=1$ and $n^*=n$

\item[(3)] Update each ensemble member
$$u_{n+1}^{(j)}=u_{n}^{(j)}+\C_{n}^{u\cG}(\C_{n}^{\cG \cG} +\alpha_{n}\Gamma  )^{-1}(y+\sqrt{\alpha_n}\xi_{n}-\cG_{n}^{(j)}),\qquad j\in\{1,\dots, J\}$$
where 
\begin{align}\label{equ3}
C_{n}^{\cG\cG}\equiv &\frac{1}{J-1}\sum_{j=1}^{J}(\cG_{n}^{(j)}-\overline{\cG}_{n})\otimes (\cG_{n}^{(j)}-\overline{\cG}_{n})\\
C_{n}^{u\cG} \equiv & \frac{1}{J-1}\sum_{j=1}^{J} (u_{n}^{(j)}-\overline{u}_{n})\otimes (\cG_{n}^{(j)}-\overline{\cG}_{n})\label{equ3B}
\end{align}

with $\overline{u}_{n}\equiv \frac{1}{J}\sum_{j=1}^{J}u_{n}^{(j)}$ and $\overline{\cG}_{n}\equiv \frac{1}{J}\sum_{j=1}^{J}\cG_{n}^{(j)}$.
\end{itemize}
$n\gets n+1$

}
        \SetKwInOut{Output}{output}
\Output{$\{u_{n^*}^{(j)}\}_{j=1}^{J}$ converged ensemble}
\caption{{\bf Generic Ensemble Kalman Inversion (with perturbed observations)}}\label{Al1}
\end{algorithm}

\section{Literature review} \label{Liter}
We discuss some existing regularisation strategies for EKI within the inversion setting posed in terms of the unregularised least-squares formulation of (\ref{equ1}), and which leads to the classical EKI formulation in (\ref{equ2}). We highlight that there is a new alternative EKI methodology and algorithms proposed in \cite{Tong1,Tong2} that arise from regularising (\ref{equ1}) with a Tikhonov-like term. In addition, the review is focused on PDE-constrained inversion; for a review of modern Kalman methods in the context of the data assimilation framework we refer the reader to the recent work on \cite{Jana} and references therein. 

\subsection{EKI as an iterative solver for identification problems.} \label{new_lab}
The initial versions of EKI \cite{ILS13a} for generic identification problems proposed to use the classical EnKF  \cite{doi:10.1029/94JC00572} update formula\footnote{the classical EnKF consists of eq (\ref{equ2}) with $\alpha_{n}=1$ fixed throughout the iterations} as an iterative solver for (\ref{equ1}) by introducing an artificial dynamical system. For various PDE-constrained identification problems, the work of \cite{ILS13a} numerically showed that this early version of EKI approximated well the solutions of (\ref{equ1}) (with $\mathcal{S}_{0}= \text{span}\{u_{0}^{(j)}\}_{j=1}^{J}$) within the first few iterations. However, they noted the algorithm became unstable if it was allowed to iterate after the data misfit (\ref{equ1}) had reached the noise level $\delta$ defined by 
\begin{eqnarray}\label{eqB:201}
\delta =\norm{ \Gamma^{-1/2}(y-\cG(u^{\dagger})},
\end{eqnarray}
where $u^{\dagger}$ is the truth. This lack of stability led to the work of \cite{EnsembleYo} where links between EKI and the regularising Levenberg-Marquardt (LM) scheme of \cite{Hanke} were first established and used to develop a regularising version of EKI. For the LM scheme, the work of \cite{Hanke} ensures that, under certain assumptions of the forward map $\cG$, the scheme in (\ref{eqA:14}) converges to the solution of (\ref{equ1}) (as $\delta \to 0$) provided that (i) the regularisation parameter $\alpha_{n}$ satisfies
\begin{eqnarray}\label{eqB:1}
\rho \norm{ \Gamma^{-1/2}(y-\cG(m_n))}\leq \alpha_{n}\norm{ \Gamma^{1/2}(D\cG(m_{n})\C_{n}D\cG^{*}(m_n)+\alpha_{n}\Gamma) (y-\cG(m_n))},
\end{eqnarray}
where $\rho<1$ is a tuning parameter, and (ii) that the algorithm is terminated at an iteration level $n^{*}$ determined by the following \textit{discrepancy principle}
\begin{eqnarray}\label{eqB:2}
\norm{ \Gamma^{-1/2}(y-\cG(m_{n^{*}}))}\leq \tau \delta< \norm{ \Gamma^{-1/2}(y-\cG(m_{n}))}\qquad 0\leq n<n^{*}
\end{eqnarray}
where  $\tau$ is another tuning parameter that must satisfy $\tau >1/\rho$. In  \cite{EnsembleYo}, these regularisation strategies from LM were adapted for the selection of $\alpha_{n}$ in EKI via using derivative-free Gaussian approximations in (\ref{eqB:1})-(\ref{eqB:2}). We refer to the approach from \cite{EnsembleYo} as \textbf{EKI-LM} (see Algorithm \ref{Al4}).

The numerical results of \cite{EnsembleYo} showed that \textbf{EKI-LM} enabled stability and accuracy for sufficiently large ensembles. Further work that has explored \textbf{EKI-LM} can be found in \cite{pub:28949,Neil,pub:29192} as well as some practical applications including seismic tomography \cite{seismic}, modeling of intracraneal pressure \cite{Chido} and time fractional diffusion inverse problems \cite{cited1}.

Despite of addressing stability in EKI, the approach \textbf{EKI-LM} suffers from a potential practical limitation that arises from the fact that it relies on the tuning parameters $\rho$ and $\tau$ in (\ref{EQ1}). Larger $\rho$'s yield larger $\alpha_n$ and in turn more iterations to converge. Smaller $\rho$'s, while desirable for computational efficiency, result in larger $\tau$'s which can lead to larger data misfit and possible loss of accuracy from stopping too early (via (\ref{EQ2})). Selecting these parameters in a computationally optimal fashion becomes more crucial when EKI is combined with complex parameterisations of the unknown. As we discuss in subsection \ref{parame}, there are cases for which, instead of using EKI directly on the physical property that we wish to infer, we need to parameterise the unknown to be able to capture properties that are not necessarily encoded in the initial ensemble. For example, in \cite{EnsembleYo} the LM approach for EKI was applied with a level-set parameterisation of the unknown in order to characterise discontinuous properties (i.e. discontinuous conductivity in the context of EIT). In comparison to the simpler case in which EKI directly estimates a physical property of interest,  \cite{EnsembleYo} found that not only they needed a larger ensemble to achieve converge, but also more EKI iterations. The application of \textbf{EKI-LM} with level-set parameterisations for seismic imaging in \cite{seismic} reported up to 40 iterations to achieve convergence. The numerical results of  \cite{Neil} also show that when EKI is combined with various others parameterisations of the unknown, the number of iterations of EKI can become large even for simple 1D and 2D forward modelling settings. We aim at addressing these very same issues with the selection of $\alpha_{n}$ that we propose in (\ref{eq: data-misfit controller}) and that we motivate in Section \ref{Moti} using the Bayesian perspective of EKI discussed earlier.

\begin{algorithm}[h!]

\SetAlgoLined

\SetKwInOut{Input}{Input}
\Input{Same from Algorithm \ref{Al1}, $\rho<1$, and $\tau>\frac{1}{\rho}$ as well as the noise level $\delta$.
}
Set $\{u^{(j)}\}_{j=1}^{J}=\{u_{0}^{(j)}\}_{j=1}^{J}$ and $\theta=0$\\
 \While{$\theta<1$ }{
\begin{itemize}
\item[(1)] Compute $\cG^{(j)}\equiv \cG(u^{(j)}),\qquad j\in\{1,\dots, J\}$
\item[(2)] Compute $\alpha_{n}$ such that 
\begin{eqnarray}\label{EQ1}
\rho \norm{ \Gamma^{-1/2}(y-\overline{\cG}_{n})}\leq \alpha_{n}\norm{ \Gamma^{1/2}(\C_{n}^{\cG \cG}+\alpha_{n}\Gamma)^{-1}(y-\overline{\cG}_{n})}
\end{eqnarray}
 \textbf{if} 
\begin{align}\label{EQ2} 
\norm{ \Gamma^{-1/2}(y-\overline{\cG}_{n})}\leq \tau \delta &.
\end{align}
 \textbf{then} \newline
\vspace{-2.5mm} 
\hspace{10mm}	 set $n^*=n.$\newline
\vspace{-2.5mm}

\hspace{10mm} \textbf{break}

\item[(3)] Update each ensemble member using (\ref{equ2}) (see also Step 3 from Algorithm \ref{Al1})
\end{itemize}
}
        \SetKwInOut{Output}{output}
\Output{$\{u_{n^*}^{(j)}\}_{j=1}^{J}$ converged ensemble}
\caption{{\bf (EKI-LM) EKI with LM selection of $\alpha_{n}$ and stopping criteria }}\label{Al4}
\end{algorithm}


\subsection{EKI as Gaussian approximation in linearised Bayesian tempering.}

Although the Bayesian perspective of EKI as a Gaussian approximation within the tempering setting was initially mentioned in \cite{SS17} and further developed in \cite{Iglesias_2018}, the early work of  \cite{EME2,EMERICK20133} established a strong link between EKI and the Bayesian setting which led exactly to the same EKI scheme from Algorithm \ref{Al1}. However, instead of using tempering, they used a heuristic approach in which the data was inverted multiple times with noise inflated by $\sqrt{\alpha_n}$. From algebraic manipulations they figured out that, in order to accurately sample from the posterior in the linear-Gaussian case, their inflation parameter $\alpha_{n}$ must satisfy \eqref{eqA:9} with $N=n^*$, where  $n^*$ is the total number of EKI iterations.  Note that, in the tempering setting, (\ref{eqA:9}) is simply a consequence of the definition of $\alpha_{n}$ in (\ref{eqA:8}) as well as the definitions of $t_{0}$ and $t_{N+1}$ in (\ref{eqA:6B}). Moreover,  in the general Bayesian tempering settings, expression (\ref{eqA:9}) holds for the general nonlinear case and without any assumptions on the prior. Nonetheless, those considerations for the linear-Gaussian case from \cite{EMERICK20133} enabled them to (i) justify the use of (\ref{eqA:9}) in the nonlinear case and (ii) to propose a simple selection of $\alpha_{n}$ given by $\alpha_{n}=n^{*}$ with $n^{*}$ selected a priori (which trivially satisfies (\ref{eqA:9})). This approach, referred to as Ensemble Smoother with Multiple Data Assimilation (ES-MDA), has been popularised in petroleum engineering applications (see \cite{evensen_crap} and references therein). However, from the insight gained from the link between EKI and the LM scheme, this selection of $\alpha_{n}$ was discouraged in \cite{Iglesias2014} since the stability of EKI, as shown in \cite{Iglesias2014, EnsembleYo}, requires $\alpha_{n}$ to be large at the beginning of the iterations, and gradually decrease as the data misfit approaches the noise level. Further versions of ES-MDA \cite{other} adopted selections of $\alpha_{n}$ similar to those initially proposed in \cite{EnsembleYo} based on the LM scheme.

As discussed in the previous section, a selection of $\alpha_{n}$ based on the adaptive-tempering Sequential Monte Carlo (SMC) method of \cite{Kantas} is proposed in \cite{Iglesias_2018}. More specifically, $\alpha_{n}$ is selected so that
\begin{eqnarray}\label{eq36}
\Bigg[\sum_{j=1}^{J}(\mathcal{W}_{n}^{(j)}[\alpha_n])^2\Bigg]^{-1}= J_{*},
\end{eqnarray}
where $J_{*}$ is a tuning user-defined parameter and 
\begin{eqnarray}\label{eq36B}
\mathcal{W}_{n}^{(j)}[\alpha_n]= \frac{\exp\Big[-\frac{1}{2}\norm{ (\alpha_{n}\Gamma)^{-1/2}(y-\cG(u_{n}^{(j)}))}^2\Big]}{\sum_{s=1}^{J} \exp\Big[-\frac{1}{2}\norm{ (\alpha_{n}\Gamma)^{-1/2}(y-\cG(u_{n}^{(s)}))}^2\Big] }.
\end{eqnarray}
The left hand side of (\ref{eq36}) is the ESS of the ensemble approximation of $\mu_{t_n}$. For further details we refer the reader to \cite{Iglesias_2018}, where the selection of $\alpha_{n}$ according to (\ref{eq36}) was implemented in a batch-sequential EKI framework to sequentially solve an inverse problem that arises in resin transfer moulding.  The same EKI methodology was applied  in \cite{DESIMON2018220} for parameter identification of the heat equation, including identification of thermal conductivity and heat capacitance given boundary measurement of heat flux. Both time-dependent applications tackled in \cite{DESIMON2018220,Iglesias_2018} involved the inversion of small number of measurements (e.g. $<30$) at each observation time, and with relatively large noise informed by measurement protocols. However, as discussed in the previous section, the selection of $\alpha_{n}$ via solving (\ref{eq36}) becomes problematic for large number of observations and/or small observational noise, since (\ref{eq36B}) requires the computation of the likelihood (\ref{eqA:7}). Another limitation is that the efficiency of the approach used in \cite{DESIMON2018220,Iglesias_2018} relies on the tuning parameter $J_{*}$ in (\ref{eq36}).



\subsubsection{EKI as a discretisation of Stochastic Differential Equations}

The pioneering work of \cite{SS17} has shown that EKI can be derived as a discretisation of an SDE system for the ensemble of particles in EKI. This so-called continuos-time limit or seamless formulation of EKI has led to further theoretical understanding of the EKI framework \cite{SS18,Bl_mker_2019,doi:10.1137/17M1132367}. In addition, using alternative discretisation schemes of the seamless formulation of EKI, leads to new EKI algorithms in both the context of optimisation \cite{SS18,ML} and sampling within the Bayesian approach \cite{Alfredo1}. 

In the seamless formulation of EKI, the regularisation parameter $\alpha_{n}$ from the classical EKI becomes the inverse of the mesh-size/discretisation step. The work of \cite{ML} proposes to choose this parameter according to
\begin{eqnarray}\label{seam}
\alpha_{n}^{-1}=\frac{\alpha_{0}^{-1}}{\norm{U}_{F}+\epsilon}
\end{eqnarray}
where $\alpha_{0}$ and $\epsilon$ are user defined parameters, $\norm{\cdot}_{F}$ denotes the Frobenius norm and $U$ is a matrix with entries $U_{j,k}=(\cG(u_{n}^{(k)})-\overline{\cG}_{n})^{T}\Gamma^{-1}(\cG(u_{n}^{(j)})- y)$. To the best of our knowledge, this selection of $\alpha_{n}$ has only been used for new different EKI algorithms including the ones arising from Forward-Euler  \cite{ML} and the implicit split-step method \cite{Alfredo1}. While (\ref{seam}) is also a perfectly valid choice of $\alpha_{n}$ for classical EKI, we emphasise that it relies on the choice of the tuning parameters $\alpha_{0}$ and $\epsilon$.






\section{The proposed regularisation framework for EKI}\label{Moti}

The aim of this section is to motivate the new approach (data misfit controller) that we propose via \eqref{eq: data-misfit controller} to select $\alpha_{n}^{-1}$ for the classical EKI setting given in (\ref{equ2}). We motivate our approach using tempering within the Bayesian setting for inverse problems in which $\alpha_{n}^{-1}$ is the step size between consecutive tempering measures. After introducing some  notation and asumptions in subsection \ref{notation}, in subsection \ref{subsec: information gain} we define tempering in the continuous setting which we require in order to compute an approximation of Jefreys' divergence between two consecutive measures. Under this approximation, we further show that our selection of $\alpha_{n}^{-1}$ \eqref{eq: data-misfit controller} controls Jefreys' divergence. The control parameter is then selected following a heuristic approach based on a discrepancy principle that we introduced in subsection \ref{subsec: mean-variance pair}.

\subsection{Notation, definitions and assumptions}\label{notation}

We assume that the space of inputs, $\mathcal{H}$, is a real-valued separable Hilbert space. The Borel-sigma algebra over $\mathcal{H}$ is denoted by $\mathcal{B}(\mathcal{H})$. In this section we assume that the prior is a Gaussian measure $\mu_{0}=\mathcal{N}(m,\C):\mathcal{B}(\mathcal{H})\to [0,1]$ on $(\mathcal{H},\mathcal{B}(\mathcal{H}))$.
Given two probability measures  $\widehat{\mu}$ and $\mu$ on a measurable space $(\mathcal{Y},\Sigma)$, if $\widehat{\mu}$ is absolutely continuous with respect to $\mu$, we define the Kullback-Leibler divergence of $\widehat{\mu}$ with respect to $\mu$ by \cite{10.2307/2236703},
\begin{equation*}
\mathrm{D}_\mathrm{KL}( \widehat{\mu} || \mu ) := \int_\mathcal{Y} \log\left( \frac{ \mathrm{d} \widehat{\mu} }{\mathrm{d} \mu } \right) \, \mathrm{d} \widehat{\mu}
\end{equation*}
where $\frac{ \mathrm{d} \widehat{\mu} }{\mathrm{d} \mu }$ denotes the Radon-Nikodym derivative of $\widehat{\mu}$ with respect to $\mu$. For two equivalent probability measures $\widehat{\mu}$ and $\mu$ (i.e. $\widehat{\mu} \ll \mu$ and $\mu \ll \widehat{\mu}$) we define the symmetrized KL divergence or Jefrey's divergence:
\begin{equation}\label{KL_sym}
\mathrm{D}_{\mathrm{KL},2}( \mu, \widehat{\mu} ) := \mathrm{D}_\mathrm{KL}( \widehat{\mu} || \mu ) + \mathrm{D}_\mathrm{KL}( \mu || \widehat{\mu} )
\end{equation}
We consider the following assumption on the forward operator:
\begin{assumption} \label{ass: forward map}
Assume that the forward map $\mathcal{G}:\mathcal{H} \to \mathbb{R}^M$ satisfies the following two conditions:
\begin{enumerate}
\item For every $\epsilon>0$ there is an $\beta=\beta(\epsilon) \in \mathbb{R}$ such that, for all $u\in \mathcal{H}$,
\begin{equation*}
\left\| \Gamma^{-1/2}\mathcal{G}(u) \right\| \leq \exp( \epsilon \|u\|_\mathcal{H}^2 + \beta )
\end{equation*}
\item For every $r>0$ there is an $K=K(r) >0 $ such that, 
\begin{equation*}
\left\| \Gamma^{-1/2}( \mathcal{G}(u_1) - \mathcal{G}(u_2) ) \right\| \leq K \|u_1-u_2\|_\mathcal{H}
\end{equation*}
for all $u_1,u_2 \in \mathcal{H}$ with $\max \left\lbrace \|u_1\|_\mathcal{H},\,\|u_2\|_\mathcal{H}  \right\rbrace < r$.
\end{enumerate}
\end{assumption}
This assumption was used in \cite[Assumption 2.7 ]{Andrew} to establish well-posedness of the infinite-dimensional Bayesian framework. Let us recall the least-square functional $\Phi(u;y)$ defined in \eqref{equ1B}. For any $y\in \mathbb{R}^{M}$, Fernique's theorem and condition 1 in Assumption \ref{ass: forward map} implies that $\Phi(\cdot; y)^m$ is square integrable, for any $m\geq 1$, with respect to the prior $\mu_0=\mathcal{N}(m,\mathcal{C})$. Integrability of $\Phi(\cdot; y)^m$ is used through the rest of this section and for the proofs in \ref{app: appendix}. 

%




\subsection{Controlling Jeffreys' divergence between successive tempering measures} \label{subsec: information gain}

In this section we show that the data-misfit controller \eqref{eq: data-misfit controller} ensures that the level of information between any two consecutive tempering measures is approximately bounded by a user-defined threshold $\theta>0$. We measure the level of information between these measures via Jefreys' divergence defined in (\ref{KL_sym}). Our aim is to show that the selection of $\alpha_{n}^{-1}=t_{n+1}-t_{n}$ via \eqref{eq: data-misfit controller} yields 
\begin{eqnarray}\label{new1}
\mathrm{D}_{\mathrm{KL},2}(\mu_{t_{n+1}},\mu_{t_n}) \leq \theta, \qquad n\in\{1,\dots,N\}
\end{eqnarray}
In subsection \ref{subsec: mean-variance pair} we give heuristic motivations for the choice of $\theta$ based on a statistical discrepancy principle.

Our first step is the explicit characterisation of $\mathrm{D}_{\mathrm{KL},2}(\mu_{t_{n+1}},\mu_{t_n})$ which we latter approximate in order to propose our adaptive selection of $\alpha_{n}^{-1}$. This approximation, require us to work on the continuous tempering setting which we now introduce. 
\subsubsection{Continuous tempering.} For $t\in [0,1]$, we define measure $\mu_{t}$ such that the Radon-Nikodym derivative of $\mu_t$ with respect to the prior $\mu_0=\mathcal{N}(m,\C)$ satisfies, for almost every $u \in \mathcal{H}$,
\begin{equation} \label{eq: Radon-Nikodym derivative mu_t mu_0}
\frac{\mathrm{d}\mu_t}{\mathrm{d}\mu_0}(u) = N_t^{-1}\exp\left( - t\Phi(u;y)  \right)
\end{equation}
where $N_t$ is the normalizing constant defined by 
\begin{equation} \label{eq: define N_t}
N_t = \int_\mathcal{H} \exp(-t\Phi(u;y)) \, \mu_0(\mathrm{d}u)
\end{equation}
The following Lemma establishes the well-posedness of $\{\mu_{t}\, \vert\, t\in [0,1]\}$. The proof can be found in \ref{proof_lema}.
\begin{lemma} \label{lem: equivalence of TPM}
The probability measure $\mu_t$ from \eqref{eq: Radon-Nikodym derivative mu_t mu_0} exists on the probability space $(\mathcal{H}, \mathcal{B}(\mathcal{H}), \mu_0 )$ and is equivalent to $\mu_0$, i.e. $\mu_t \ll \mu_0$ and $\mu_0 \ll \mu_t$.
\end{lemma}
Using integrability of the data misfit with respect to the prior, and the equivalence proven in the previous lemma, let us define the mean and variance of the data misfit:
\begin{equation}\label{eq:mean_var}
\left\langle \Phi \right\rangle_t := \int_\mathcal{H} \Phi(w;y)\, \mu_t(\mathrm{d}w) \qquad \left\langle \Phi,\Phi \right\rangle_t := \int_\mathcal{H} \left( \Phi(w;y) - \left\langle \Phi \right\rangle_t \right)^2 \, \mu_t(\mathrm{d}w)
\end{equation}
for all $t\in [0,1]$. 

It is not difficult to see that, given $0=t_0<t_1< \cdots < t_{N+1}= 1$, Lemma \ref{lem: equivalence of TPM} also establishes the equivalence of the probability measures $\left\lbrace \mu_{t_n} : n=0,1,...,N+1 \right\rbrace $. Furthermore, from the integrability of the data misfit, it follows that the Jeffreys' divergence between consecutive measures is well-defined and, using its definition in \eqref{KL_sym}, it can be written as
\begin{equation} \label{eq: D_KL_2 explicit}
\mathrm{D}_{\mathrm{KL},2}( \mu_{t_{n}} , \mu_{t_{n+1}} ) = \alpha_{n}^{-1}\left\langle \Phi \right\rangle_{t_n}-\alpha_{n}^{-1}\left\langle \Phi \right\rangle_{t_{n+1}}.
\end{equation}
The finite-dimensional version of \eqref{eq: D_KL_2 explicit} was proven in \cite[Proposition 3.2]{ps} in the context of path sampling. 

Let us note that, since $\Phi\ge 0$, 
\begin{equation} \label{eq: D_KL_2 explicit2}
\mathrm{D}_{\mathrm{KL},2}( \mu_{t_{n}} , \mu_{t_{n+1}} ) \leq \alpha_{n}^{-1}\left\langle \Phi \right\rangle_{t_n} 
\end{equation}

\subsubsection{Approximating $\mathrm{D}_{\mathrm{KL},2}( \mu_{t_{n}} , \mu_{t_{n+1}} )$}

In practice, the expression for $\mathrm{D}_{\mathrm{KL},2}( \mu_{t_{n}} , \mu_{t_{n+1}} )$ provided in \eqref{eq: D_KL_2 explicit} cannot be used in (\ref{new1}) to find $\alpha_{n}^{-1}$ since, at time $t_{n}$,  measure $\mu_{t_{n+1}} $ is unknown. We note that if $ \left\langle \Phi \right\rangle_{t_{n+1}} \ll \left\langle \Phi \right\rangle_{t_{n}} $, then the term $\left\langle \Phi \right\rangle_{t_n} \geq 0$ in the right-hand side of \eqref{eq: D_KL_2 explicit} can be neglected and hence 
\begin{equation} \label{eq: D_KL_2 explicit3}
\mathrm{D}_{\mathrm{KL},2}( \mu_{t_{n}} , \mu_{t_{n+1}}) \approx  \alpha_{n}^{-1}\left\langle \Phi \right\rangle_{t_n}\leq \theta \Rightarrow  \alpha_{n}^{-1}\leq \frac{\theta}{\left\langle \Phi \right\rangle_{t_n}}
\end{equation}
However, approximating $\mathrm{D}_{\mathrm{KL},2}( \mu_{t_{n}} , \mu_{t_{n+1}})$ with its upper bound in \eqref{eq: D_KL_2 explicit2} may not be accurate when  $ \left\langle \Phi \right\rangle_{t_{n+1}} \approx \left\langle \Phi \right\rangle_{t_{n}} $ which is likely to happen when $\alpha_{n}^{-1}$ is small. We address this case via a first order approximation of $\left\langle \Phi \right\rangle_t$ (as a function of $t$), around $t = t_{n}$.  We first need the following proposition that we prove in \ref{proof_of_theo} and which extends, to the present infinite-dimensional setting, the results from \cite[Section 2]{tunning} for tuning tempered transitions. 

\begin{theorem} \label{the: dynamic equation}
For any bounded $t \geq 0$, let $u_t \sim \mu_t$ be an $\mathcal{H}$-valued random variable, where $\mu_t$ is the probability measure determined via formula \eqref{eq: Radon-Nikodym derivative mu_t mu_0}. For any square integrable functional $g:\mathcal{H} \to \mathbb{R}$, the expected value $\mathbb{E}\left\lbrace g(u_t) \right\rbrace$, as a function of $t$, is differentiable and its derivative $\mathbb{E}\left\lbrace g(u_t) \right\rbrace'$ satisfies:
\begin{equation} \label{eq: dynamic equation}
\mathbb{E}\left\lbrace g(u_t) \right\rbrace'=-\mathrm{Cov}\left\lbrace g(u_t),\Phi(u_t;y) \right\rbrace
\end{equation}
where $\mathrm{Cov}$ denotes covariance between two random variables. 
\end{theorem}
The following corollary is a simple consequence of Theorem \ref{the: dynamic equation} with $g(\cdot)=\Phi(\cdot, y)$.
\begin{corollary} \label{pro: energy fluctuation}
For any $t \in [0,1]$, let $\mu_t$ denote the probability measure determined via the tempering setting \eqref{eq: Radon-Nikodym derivative mu_t mu_0}. The mean $\left\langle \Phi \right\rangle_t$ defined in (\ref{eq:mean_var}), as a function of $t\in[0,1]$, is differentiable. Moreover, its derivative is given by 
\begin{equation*}
\left\langle \Phi \right\rangle_t' = -\left\langle \Phi,\Phi \right\rangle_t 
\end{equation*}
\end{corollary}

We now come back to the case in which $\alpha_{n}^{-1}= t_{n+1} -t_{n}$ needs to be small so that \eqref{new1} can be satisfied for a given threshold $\theta>0$. Using Proposition \ref{pro: energy fluctuation} we have
$$ \left\langle \Phi \right\rangle_{t_{n+1}} \approx \left\langle \Phi \right\rangle_{t_{n}} + \alpha_{n}^{-1}\left\langle \Phi \right\rangle_{t_{n}}' = \left\langle \Phi \right\rangle_{t_{n}} - \alpha_{n}^{-1}\left\langle \Phi,\Phi \right\rangle_{t_{n}}. $$
Thus, 
\begin{equation}\label{new2}
\mathrm{D}_{\mathrm{KL},2}( \mu_{t_{n}} , \mu_{t_{n+1}}) \approx  \alpha_{n}^{-2} \left\langle \Phi,\Phi \right\rangle_{t_{n}} \leq \theta  \Rightarrow  \alpha_{n}^{-1}\leq \sqrt{\frac{\theta}{\left\langle \Phi,\Phi \right\rangle_{t_{n}} }}.
\end{equation}
We postulate that the approximation used in (\ref{new2}) is reasonable provided that it honours inequality \eqref{eq: D_KL_2 explicit2}, i.e. 
$$\mathrm{D}_{\mathrm{KL},2}( \mu_{t_{n}} , \mu_{t_{n+1}}) \approx  \alpha_{n}^{-2} \left\langle \Phi,\Phi \right\rangle_{t_{n}} \leq \alpha_{n}^{-1}\left\langle \Phi \right\rangle_{t_n} $$
which in turn gives us the  following condition on $\alpha_{n}^{-1}$
$$\alpha_{n}^{-1}\leq \frac{\left\langle \Phi \right\rangle_{t_n} }{\left\langle \Phi,\Phi \right\rangle_{t_{n}}}$$
Whenever this condition is not satisfied, we approximate $\mathrm{D}_{\mathrm{KL},2}( \mu_{t_{n}} , \mu_{t_{n+1}})$ via its upper bound as in (\ref{eq: D_KL_2 explicit3}). In summary we propose the following approximation 
\begin{equation*}
\mathrm{D}_{\mathrm{KL},2}( \mu_{t_{n}} , \mu_{t_{n+1}}) \approx \begin{cases} 
\alpha_{n}^{-1}\left\langle \Phi \right\rangle_{t_{n}} & \text{if} \, \alpha_{n}^{-1} \geq \frac{\left\langle \Phi \right\rangle_{t_{n}}}{ \left\langle \Phi,\Phi \right\rangle_{t_{n}} } \\
\alpha_{n}^{-2}\left\langle \Phi, \Phi \right\rangle_{t_{n}}    & \text{if } \, \alpha_{n}^{-1} \leq \frac{\left\langle \Phi \right\rangle_{t_{n}}}{ \left\langle \Phi,\Phi \right\rangle_{t_{n}} }
  \end{cases}
\end{equation*}
which can be written as follows,
\begin{equation} \label{eq: D_KL_2 approximate}
\mathrm{D}_{\mathrm{KL},2}( \mu_{t_{n}} , \mu_{t_{n+1}}) \approx \min \left\lbrace \alpha_{n}^{-1}\left\langle \Phi \right\rangle_{t_{n}}, \, \alpha_{n}^{-2}  \left\langle \Phi, \Phi \right\rangle_{t_{n}} \right\rbrace
\end{equation}
It is not difficult to see that if we define $\alpha_{n}^{-1}$:
\begin{equation} \label{eq: data-misfit controller1}
\alpha_{n}^{-1}= \min \left\lbrace \max \left\lbrace \frac{\theta}{\left\langle \Phi \right\rangle_{t_{n}}} ,\, \sqrt{ \frac{\theta}{\left\langle \Phi, \Phi \right\rangle_{t_{n}} } }\right\rbrace,\, 1 -t_{n}  \right\rbrace,
\end{equation}
then \eqref{new1} is satisfied using the approximation \eqref{eq: D_KL_2 approximate}. In addition, let $n^{*}$ be the first iteration such that 
$$\alpha_{n^{*}}^{-1} = 1 -t_{n^*}  =1-t_{n^*-1}+\alpha_{n^*-1}^{-1}=1-\sum_{j=1}^{n^*-1}\alpha_{j}^{-1}$$
We note that \eqref{eqA:9} is satisfied with $n^*=N$; hence $n^{*}$ defines the natural stopping iteration level.
 
Expression \eqref{eq: data-misfit controller1} has the same form of the data misfit controller in \eqref{eq: data-misfit controller} with (i) $\left\langle \Phi \right\rangle_{t_{n}}$ and $\left\langle \Phi, \Phi \right\rangle_{t_{n}}$ approximated with empirical mean and variance from the ensemble of particles (recall $u_{n}^{(j)}$ are approximate samples of $\mu_{t_{n}})$, and (ii) the choice $\theta = M/2$ which we motivate from heuristic arguments in the next subsection. 


\subsection{A statistical discrepancy Principle.} \label{subsec: mean-variance pair}
According to Morozov's discrepancy principle \cite{doi:10.1002/zamm.19860660714}, if $u\in \mathcal{H}$ satisfies 
$$\norm{\Gamma^{-1/2}(y-\cG(u))}\leq \delta,$$
where $\delta$ is the noise level defined in $ \eqref{eqB:201}$, then $u$ is an acceptable estimate of the classical (deterministic) inverse problem (\ref{equ0}). In the context of iterative regularisation \cite{Iterative}, a form of the discrepancy principle is often used as an early stopping rule (e.g see expression \eqref{eqB:2}) to avoid instabilities. In the Bayesian setting, however, the noise level is a random variable. In fact, under the assumption of additive centred Gaussian noise (\ref{equ0}), the squared noise level $\delta^2$ (see \eqref{eqB:201}) is a realisation of a chi-square random variable with  $M$ degrees of freedom. Therefore, its mean and variance are:
$$\mathbb{E}\left\lbrace \delta^2 \right\rbrace = M,\qquad  \mathrm{Var}\left\lbrace \delta^2 \right\rbrace = 2M.$$
Inspired by classical discrepancy principle, we then postulate that an acceptable estimate, $u$, must satisfy either one of the following conditions
\begin{description}
\item[C1] (accuracy) 
\begin{equation} \label{eq: accuracy test}
\mathbb{E}\left\lbrace \left\| \Gamma^{-1/2}\left( y - \mathcal{G}\left( u \right) \right) \right\|^2 \right\rbrace \leq \mathbb{E}\left\lbrace \delta^2 \right\rbrace \equiv M. 
\end{equation}
\item [C2] (uncertainty)
\begin{equation} \label{eq: uncertainty test}
\mathrm{Var}\left\lbrace \left\| \Gamma^{-1/2}\left( y - \mathcal{G}\left(u \right) \right) \right\|^2 \right\rbrace  \leq \mathrm{Var}\left\lbrace \delta^2 \right\rbrace \equiv 2M.
\end{equation}
\end{description}
We invoke this statistical discrepancy principle to determine the step size $\alpha_n^{-1}$ of the $n$th Bayesian sub-problem \eqref{eqA:9B} that arises from transition between tempering measures $\mu_{t_{n}}$ and $\mu_{t_{n}}$. More specifically, let us consider the following remark.
\begin{remark}\label{rema0}
We can view (\ref{eqA:7}) as an iterative application of Bayes rule, where at the $n$th iteration level we have a prior $\mu_{t_n}(u)$ and a likelihood defined by the observational model 
\begin{eqnarray}\label{eqA:9B}
y=\cG(u)+\sqrt{\alpha_{n}}\eta,\qquad \eta \sim N(0,\Gamma).
\end{eqnarray}
In other words, (\ref{eqA:7})  defines a sequence of Bayesian inverse problems similar to the original one\footnote{recall the observational model for original problem is $y=\cG(u)+\eta$ with $\eta\sim N(0,\Gamma)$).} but with a Gaussian error $\sqrt{\alpha_n}\eta$ that has a covariance matrix $\Gamma$ inflated by $\alpha_{n}$. We can then think of $\mu_{t_{n+1}}$ given by (\ref{eqA:7}) as the distribution of $u\vert y$ under the observational model (\ref{eqA:9B}) and prior $\mu_{t_{n}}$. 
\end{remark}
We apply \eqref{eq: accuracy test} and \eqref{eq: uncertainty test} to the sub-problem in Remark \ref{rema0}. To this end, we use $\alpha_{n} \Gamma$  and $u_{n}\sim \mu_{t_{n}}$ instead of $\Gamma$ and $u$ in formulas \eqref{eq: accuracy test} and \eqref{eq: uncertainty test} to obtain
\begin{equation} \label{eq: accuracy condition}
2\alpha_{n}^{-1} \left\langle \Phi \right\rangle_{t_{n}} =  \mathbb{E} \left\lbrace \left\| (\alpha_{n}\Gamma)^{-1/2}\left( y  - \mathcal{G}( u_n ) \right) \right\|^2 \right\rbrace \leq M
\end{equation}
\begin{equation} \label{eq: uncertainty condition}
4\alpha_{n}^{-2}  \left\langle \Phi, \Phi \right\rangle_{t_{n}}  = \mathrm{Var} \left\lbrace \left\| (\alpha_{n}\Gamma)^{-1/2}\left( y  - \mathcal{G}( u_n ) \right) \right\|^2 \right\rbrace  \leq 2M
\end{equation}
By choosing $\alpha_{n}^{-1}$ via (\ref{eq: data-misfit controller1}) with $\theta=M/2$, we enforce that at least one of the above condition \eqref{eq: accuracy condition} or \eqref{eq: uncertainty condition} are satisfied, while ensuring that  $t_{n+1} = t_{n}+\alpha_{n}^{-1}\leq 1$ as required. In other words, the proposed statistical discrepancy principle applied to each intermediate Bayesian sub-problem arising from tempering, yields the data misfit controller \eqref{eq: data-misfit controller} which, from the previous subsection, controls Jefreys' divergence with threshold $\theta=M/2$. We incorporate the data-misfit controller for the selection of $\alpha_{n}^{-1}$ in EKI; we summarise the scheme in Algorithm \ref{Al2} ({\bf EKI-DMC}). 
\begin{remark}
For the applications that we discuss in Section \ref{Num_testing}, we have found that the mean of the data misfit $\overline{\Phi}_{n}$ is considerably smaller than the standard deviation $\sigma_{\Phi_n}$. Therefore, 
\begin{equation} \label{ss}
 \frac{\overline{\Phi}_{n}^2}{\sigma_{\Phi_n}^2}\leq 1< \frac{M}{2}, \Rightarrow \frac{M}{ 2 \overline{\Phi}_{n}}>   \sqrt{ \frac{M}{2\sigma_{\Phi_{n}}^2 }} 
 \end{equation}
Hence, in this case, $\alpha_{n}^{-1}$ is given by condition (\ref{eq: accuracy condition}) rather than condition \eqref{eq: uncertainty condition}. Intuitively, \eqref{ss} is likely to occur in problems where the prior is sufficiently wide so that it comprises the truth. However, for a prior that is centred too far from the truth and with small variance, we would then expect for $\overline{\Phi}_{n}$ to be larger than $\sigma_{\Phi_n}$. In this case, DMC selects $\alpha_{n}^{-1}$ according to \eqref{eq: uncertainty condition} in order to allow for larger steps in the tempering settings. 

\end{remark}

\begin{algorithm}[h!]

\SetAlgoLined

$\Phi_{n}\equiv \{\Phi(u_{n}^{(j)};y)\}_{j=1}^{J}$. More specifically, we select $\alpha_{n}^{-1}$ via

\SetKwInOut{Input}{Input}
\Input{Same from Algorithm \ref{Al1}
}
Set $n=0$ and $t_{0}=0$\\
 \While{$t_{n}<1$ }{
\begin{itemize}
\item[(1)] Compute $\cG_{n}^{(j)}= \cG(u_{n}^{(j)}),\qquad j\in\{1,\dots, J\}$
\item[(2)] Compute $\alpha_{n}$ via (\ref{eq: data-misfit controller})
 \item[(3)] Update each ensemble member using (\ref{equ2}) (see also Step 3 from Algorithm \ref{Al1})

\end{itemize}
$n\gets n+1$\\
$t_n\gets t_{n}+\alpha_n^{-1}$

}
        \SetKwInOut{Output}{output}
\Output{$\{u_{n^*}^{(j)}\}_{j=1}^{J}$ converged ensemble}
\caption{{\bf EKI with Data Misfit Controller (EKI-DMC)}}\label{Al2}
\end{algorithm}




\section{Numerical Testing}\label{Num_testing}

In this section we test the performance of {\bf EKI-DMC} in the context of Electrical Impedance Tomography (EIT) with the Complete Electrode Model (CEM) that we introduce in subsection \ref{tests}. We use two parameterisations of the unknown that we introduce in subsection \ref{parame}. Implementation aspects and measures of performance are discussed in subsections \ref{imple}-\ref{per}. In subsection \ref{eit_smooth}-\ref{eit_disc} we discuss numerical results which includes comparing the performance of {\bf EKI-DMC} and {\bf EKI-LM}.

\subsection{Complete Electrode Model (CEM)}\label{tests}

Given electric currents $\{I_{k}\}_{k=1}^{m_{e}}$ injected through a set of surface electrodes $\{e_{k}\}_{k=1}^{m_{e}}$ placed on the boundary, $\partial D$, of $D$, the CEM consist of finding $[v,\{V_{k}\}_{k=1}^{m_{e}}]$ where $v$ is the voltage in $D$ and $\{V_{k}\}_{k=1}^{m_{e}}$ are the voltages on the electrodes. The dependent variables  $[v,\{V_{k}\}_{k=1}^{m_{e}}]$ are given by the solution to \cite{Cheney}:
\begin{eqnarray}\label{cem1}
&~~~~~~~\nabla \cdot \kappa \nabla v&=0\qquad\textrm{in}~~D,\\
&v+z_{k} \kappa \nabla v\cdot \n&=V_{k} \qquad\textrm{on}~~e_k, ~~k=1,\dots,m_{e},\\
&~~~~~~~~~\nabla v\cdot \n &=0 \qquad\textrm{on}~~\partial D\setminus \cup_{k=1}^{m_{e}}e_{k},\\
&\int_{e_{k}}\kappa \nabla v\cdot \n ~ds &= I_{k}\qquad ~~k=1,\dots,m_{e},\label{cem2}
\end{eqnarray}
where $\kappa$ is the electric conductivity of $D$ and $\{z_{k}\}_{k=1}^{m_{e}}$ are the electrodes' contact impedances. We consider an experimental setting consisting of $n_{p}$ current patterns $I_{1}=\{I_{1,k}\}_{k=1}^{m_{e}},\dots I_{n_{p}}=\{I_{n_p,k}\}_{k=1}^{m_{e}}$. For each of these patters $\{I_{j,k}\}_{k=1}^{m_{e}}$, we denote by $\{V_{j,k}\}_{k=1}^{m_{e}}$ the prediction of voltages at the electrodes defined by the CEM (\ref{cem1})-(\ref{cem2}). For simplicity we assume that the contact impedances of the electrodes are known. We define the map
\[
\mathcal{F}(\kappa)=V\equiv \big[\{V_{1,k}\}_{k=1}^{m_{e}}, \dots , \{V_{n_p,k}\}_{k=1}^{m_{e}}\big],
\]
that for every conductivity, produces voltage measurements. EIT consist of finding the conductivity $\kappa^{\dagger}$ of $D$ given measurements of $V^{\dagger}=\mathcal{F}(\kappa^{\dagger})$. For a review of the EIT problem we refer the reader to \cite{EIT_revew}. 

\subsection{Parameterisations of EKI}\label{parame}

In this subsection we introduce two maps, $\mathcal{P}_{1}$ and $\mathcal{P}_{2}$, that we use to parameterise the unknown physical property $\kappa^{\dagger}$ that we wish to estimate using the EKI framework introduced earlier. 

\subsubsection{Parameterisation $\mathcal{P}_{1}$. Smooth properties}\label{smoothpara}

Let us first discuss the case in which prior knowledge suggest that the unknown physical property $\kappa^{\dagger}$ is a smooth function. For simplicity we only consider the 2D case that we use in the numerical experiments of Section \ref{EIT_num} but we emphasise that the approach can be used for 1D and 3D settings.  Let us introduce a \textit{Whittle-Matern (WM) parameterisation} of the unknown $\kappa(x)$ that we wish to infer via EKI. The WM parameterisation involves a positive smoothness parameter denoted by $\nu$, an amplitude scale $\sigma$, and two intrinsic lengthscales, $L_{1}$ and $L_{2}$, along  the horizontal and vertical direction, respectively. Given $\Theta=(\nu, \sigma, L_{1},L_{2} )\in \mathbb{R}_{+}\times \mathbb{R}\times \mathbb{R}_{+}^{2}$, we define an operator $\mathcal{W}_{\Theta}$ that maps every $\omega\in H^{-1-\epsilon}(D)$ (with $\epsilon>0$ arbitrary) to $\Psi=\mathcal{W}_{\Theta}\omega$, that satisfies the following fractional PDE in the domain $D$
\begin{align}\label{para5_1}\Big(\mathbb{I}-\nabla \cdot  \text{diag}(L_{1}^2,L_{2}^2)\nabla \Big)^{(\nu+1)/2}\Psi=4 \sigma^2 \pi\frac{\Gamma(\nu+1)}{\Gamma(\nu)}\sqrt{L_{1}L_{2}}\, \omega,
\end{align}
with Robin boundary conditions on $\partial D$:
\begin{align}\label{para5_1B}
\Psi-\zeta_{R} \,\text{diag}(L_{1}^2,L_{2}^2)\nabla \Psi\cdot \mathbf{n}=0.
\end{align}
where $\zeta_{R}$ is a tuning parameter. In (\ref{para5_1}) $\Gamma$ denotes the gamma function, $\mathbb{I}$ is the identity operator, and
$$ \text{diag}(L_{1}^2,L_{2}^2)\equiv \left(\begin{array}{cc}
L_{1}^2 & 0\\
0&  L_{2}^2\end{array}\right).$$
The proposed WM parameterisation of $\kappa$ is defined by
\begin{align}\label{para5_2}
\kappa= \mathcal{P}_{1}(\lambda, \Theta,\omega)\equiv  \lambda\exp\big(\mathcal{W}_{\Theta}\omega \big),
\end{align}
where $\lambda$ a positive scaling factor for $\kappa$. Note that (\ref{para5_2}) enforces that the electrical conductivity is positive. We can succinctly write (\ref{para5_2}) as, 
\begin{align}\label{para5_2B}
\kappa=\mathcal{P}_{1}(u),\qquad \text{where}\qquad u=(\lambda, \Theta,\omega)=(\lambda,\nu, \sigma, L_{1},L_{2},\omega).
\end{align}
%
%
%
%
%
The motivation behind this parameterisation comes from the theory of Gaussian random fields (GRFs). In particular, the work of  \cite{Lindgren} that shows that if $\omega\in H^{-1+\epsilon}$ is Gaussian white noise (i.e. $\omega\sim N(0,\mathbb{I})$) then 
$$\log{\kappa}=\log(\lambda)+\mathcal{W}_{\Theta}\omega\sim N(\log{\lambda},C_{\Theta})$$
 where $C_{\Theta}$ is the covariance operator induced by the Matern autocorrelation function defined by
\[
\text{ACF}_{\Theta}(x) = \sigma^2\frac{1}{2^{\nu-1}\Gamma(\nu)}\norm{x}_{L_{1},L_{2}}^\nu K_\nu \Big(\norm{x}_{L_{1},L_{2}}\Big),
\]
where $K_\nu$ is the modified Bessel function of the second kind of order $\nu$, and 
\[
\norm{x}_{L_{1},L_{2}}\equiv \sqrt{\frac{x_{1}^2}{L_{1}^2}+\frac{x_{2}^2}{L_{2}^2}}
\]
It can also be shown that if $\log{\kappa}\sim N(\log{\lambda},C_{\Theta})$, then almost surely $\log{\kappa}\in H^{\nu-\epsilon}(D)$ \cite{Matt} which further shows the role of the smoothness parameter $\nu$. 

Our choice for the boundary conditions (BCs) in (\ref{para5_1B}) follows from the work of \cite{whittlematern} that shows that Robin BCs are better suited, compared to Neumman and Dirichlet, to alleviate (via the appropriate choice of $\zeta_{R}$) undesirable boundary effects which arise from the discretisation of GRFs. Let us reiterate that our goal here is to introduce parameterisations for EKI and then suitable initial ensembles on the parameters. Whether the initial ensemble yields Gaussian (or log-Gaussian) properties is not our main concern. However, it would not be advisable to select parameterisations that, for example, restrict the values of the physical property near the boundary which is something that we would expect if $\zeta_{R}=0$ in (\ref{para5_1B}). 


The WM parameterisation in (\ref{para5_2}) allows us to incorporate the smoothness and lengthscales of the underlying field as part of the unknown that we can estimate with EKI. In contrast to the work of \cite{Neil} where the lengthscales are estimated under isotropic assumptions, here we consider the anisotropic case. We focus on vertical/horizontal anisotropy but a rotation matrix can be further introduced in (\ref{para5_1}) to characterise properties with an arbitrary (and unknown) direction of anisotropy \cite{whittlematern}

A straightforward approach to generate the initial ensemble for EKI with the parameterisation from (\ref{para5_2}) is to specify (hyper prior) densities, $\pi_{\lambda}$, $\pi_{L_{1}}$, $\pi_{L_{2}}$, $\pi_{\sigma}$ and $\pi_{\nu}$, and produce samples:
\begin{align}\label{para5_5}
u_{0}^{(j)}\equiv (\lambda^{(j)},\nu^{(j)}, \sigma^{(j)},L_{1}^{(j)},L_{2}^{(j)},\omega^{(j)}),\sim  \pi_{\lambda}\otimes \pi_{\nu}\otimes \pi_{\sigma}\otimes \pi_{L_{1}}\otimes \pi_{L_{2}}\otimes N(0,\mathbb{I})
\end{align}

\begin{remark}[KL expansion]
For some geometries of $D$, the eigenfunctions and eigenvalues of  a Matern covariance $C_{\Theta}$ are available in closed form \cite{Matt}. In this case, an equivalent formulation of the WM parameterisations can be defined in terms of the spectral decomposition of $C_{\Theta}$ (i.e. KL expansion under the prior). Of course, for computational purposes WM parameterisations can be defined on a larger (simple) domain that encloses $D$ and simply restrict the physical property to the domain of interest. While the spectral/KL approach is more standard, the approach we use here based on the work of  \cite{Lindgren} allows to naturally define the paramerisation in (\ref{para5_2B}) via (\ref{para5_1}). 
\end{remark}
\begin{remark}[Smoothness parameter and amplitude scales]\label{rem_smooth}
For simplicity, in this work the smoothness parameter $\nu$ and the amplitude scales $\sigma$ in the WM are fixed, i.e. we leave these parameters our of the inversion. For the experiments that we discuss in the subsequent sections, the spatial variability in the unknown can be captured quite effectively only via $\omega(x)$. However, we recognise that including a variable $\sigma$ can be beneficial in some cases as reported in \cite{whittlematern}. Similarly, in the context of the experiments reported later, where the aim is mainly to recover medium anomalies, we find that sensible (fixed) choices of $\nu$ are sufficient, provided that the lengthscales are properly estimated via EKI. Nevertheless, as shown in \cite{Neil}, EKI can be used to estimate the smoothness parameter $\nu$. 
\end{remark}

\subsubsection{Parameterisation $\mathcal{P}_{2}$. Piecewise-constant functions.}\label{case2}
In order to characterise piecewise-smooth functions we use the level-set approach initially proposed in \cite{refId0} for deterministic inverse problems and, more recently, for fully Bayesian \cite{levelBayes,Matt} and EKI \cite{EnsembleYo} settings. Our main modelling assumptions is that the unknown property has a background value potentially heterogeneous, and that possible anomalies/defects consist of regions with (also possibly heterogeneous) higher/lower values  than those in the background field. More specifically, let us first define the parameterisation 
\begin{align}\label{para5_10}
\kappa(x)= \mathcal{H}(\{\kappa_{\iota}\}_{\iota\in \{l,b,h\}},f)\equiv\left\{\begin{array}{cc}
\kappa_{l}&  f(x)\leq \zeta_{1}\\
\kappa_{b}& \zeta_{1}<f(x)\leq \zeta_{2}\\
\kappa_{h} & f(x)>\zeta_{2}
\end{array}\right.
\end{align}
where $\kappa_{b}$, $\kappa_{l}$ and $\kappa_{h}$ denote the background, low-value and high-value fields that characterise the unknown physical property. For simplicity we assume these variables take only constant values. The level-set function denoted by $f(x)$ determines the background $\Omega_{b}\equiv  \{x:\zeta_{1}<f(x)\leq \zeta_{2}\}$ as well as the region of low $\Omega_{l}\equiv  \{x:f(x)\leq \zeta_{1}\}$ and high $\Omega_{h}\equiv \{x\in f(x)>\zeta_{2}\}$ values. We assume that $f$ is a smooth fields with some variability that we enforce via a second level of parameterisation in terms of $\mathcal{P}_1$ introduced earlier. More specifically, we consider
\begin{eqnarray}\label{para5_11B}
f=\log{\lambda_{f}}+\mathcal{W}_{\Theta_{f}}\omega_{f} =\log (\mathcal{P}_{1}(u_{f})) 
\end{eqnarray}
where $\mathcal{P}_{1}$ is defined according to (\ref{para5_2}) and 
$$u_{f}=(\lambda_{f},\nu_{f},\sigma_{f}, L_{1,f},L_{2,f},\omega_{f}).$$

Combining (\ref{para5_10}) with (\ref{para5_11B}) we can write
\begin{align}\label{para5_13}
\kappa=\mathcal{P}_{2}(u),\qquad \text{with}\qquad u\equiv \Big\{\{\kappa_{\iota}\}_{\iota\in \{l,b,h\}},u_{f}\Big\},
\end{align}
where
\begin{align}\label{para5_1234}
\mathcal{P}_{2}(u) =\mathcal{H}\Big(\{\kappa_{\iota}\}_{\iota\in \{l,b,h\}},\log (\mathcal{P}_{1}(u_{f})) \Big)= \mathcal{H}\Big(\{\kappa_{\iota}\}_{\iota\in \{l,b,h\}},\log \Big(\mathcal{P}_{1}( \lambda_{f},\nu_{f},\sigma_{f},L_{1,f}, L_{2,f},\omega_{f})  \Big) \Big).
\end{align}
The selection of the initial ensemble for $u_{f}$ can be done similarly to the one in (\ref{para5_5}). This parameterisation can be extended for the case in which $\{\kappa_{\iota}\}_{\iota\in \{l,b,h\}}$ are unknown functions via using $\mathcal{P}_{1}$ to parameterise each of these functions. 


\subsection{Implementation}\label{imple}

For the numerical experiments discussed in the following subsections, Algorithms \ref{Al4} - \ref{Al2} are implemented in MATLAB. Let us recall that for these algorithms we need to construct the forward map $\cG=\mathcal{F}\circ \mathcal{P}$ where $\mathcal{F}$ is the operator induced by the CEM and $\mathcal{P}$ is any of the parameterisations defined earlier. The numerical implementation of the CEM from subsection \ref{tests} is conducted using MATLAB software EIDORS \cite{EIDORS}. The experimental setting consists of (i) a circular domain of unit radius centred at the origin, (ii) 16 surface electrodes with contact impedances of values 0.01 Ohms, (iii) an adjacent injection pattern with an electric current of 0.1 Amps, and (iv) measurements at each electrode. All these parameters are assumed known and fixed for the inversions of this section. Synthetic data are generated using the mesh from Figure \ref{Fig1} (right) with 9216 elements, while a coarser mesh of 7744 elements is used for inversions. The total number of measurements is $M=16^2=256$.

For the evaluation of the WM parameterisation $\kappa=\mathcal{P}_{1}(u)$, we solve (\ref{para5_1})-(\ref{para5_1B}) using the techniques from \cite{Lindgren} and which restrict us to the cases in which $\nu\in \mathbb{N}$. The discretisation of the operator $\mathbb{I}-\nabla \cdot  \text{diag}(L_{1},L_{2})\nabla$ with BCs from (\ref{para5_1B}) is performed via cell-centred finite differences. The PDE in (\ref{para5_1})-(\ref{para5_1B}) is solve in a square domain, equal to, or enclosing $D$ (the domain of definition for the PDE encoded in $\mathcal{F}$). The actual field $\kappa(x)$ that we pass into the CEM is an interpolation of $\mathcal{P}(u)$ on $D$. Parameterisations $\mathcal{P}_{2}$ is based on the truncation of the level-set function, $f(x)$ so the implementation is straightforward once all the fields $\kappa_{\iota}$ ($\iota\in {l,b,h}$) and $f(x)$ are computed. Given these construction of $\mathcal{G}$, the rest of the steps in Algorithms \ref{Al4} - \ref{Al2} are computed in a straightforward manner.

\subsection{Measures of performance for EKI}\label{per}

Given the ensemble $\{u_{n}^{(j)}\}$ computed via EKI at the $n$ iteration of the scheme ($n=0$ corresponds to the prior ensemble), our estimate of the unknown property $\kappa^{\dagger}$ is given by
\begin{align}\label{eq:2003A}
\kappa_{n}\equiv \mathcal{P}(\overline{u}_{n}) =  \mathcal{P}\Bigg(\frac{1}{J}\sum_{j=1}^{J}u_{n}^{(j)} \Bigg).
\end{align}
We measure the accuracy in terms of the relative error with respect to (w.r.t) the truth defined by
\begin{align}\label{eq:2003}
\mathcal{E}_{n}=\frac{\norm{\kappa_{n}-\kappa^{\dagger}}_{L^{2}(D)} }{\norm{\kappa^{\dagger}}_{L^{2}(D)} }.
\end{align}
We often visualise some transformed ensemble members (mainly for the initial ensemble $n=0$), i.e. 
\begin{align}\label{eq:2003B}
\kappa_{n}^{(j)}\equiv \mathcal{P}(u_{n}^{(j)}),\qquad j\in \{1,\dots, J\},
\end{align}
but note that our estimate $\kappa_{n}$ in (\ref{eq:2003A}) does not involve taking the average of the particles in (\ref{eq:2003B}); this would be particularly detrimental for $\mathcal{P}_{2}$ and $\mathcal{P}_{3}$ since averaging (\ref{eq:2003B}) will not preserve discontinuities. For these two parameterisations we also visualise an estimate of the level-set function given by 
\begin{align}\label{eq:2003LF}
f_{n}\equiv \log\mathcal{P}_{1}(\overline{u}_{f,n}) 
\end{align}

We additionally monitor the following data-misfit quantities
\begin{align}\label{eq:2004}
\mathcal{DM}_{1,n}=&\norm{\Gamma^{-1/2}\Big[y-\frac{1}{J}\sum_{j=1}^{J}\cG(u_{n}^{(j)})\Big]}\\
\mathcal{DM}_{2,n}=&\norm{\Gamma^{-1/2}(y-\cG(\overline{u}_{n}))}\label{eq:2005}\\
\mathcal{DM}_{3,n}=&\Bigg[\frac{1}{J}\sum_{j=1}^{J}\norm{\Gamma^{-1/2}(y-\cG(u_{n}^{(j)}))}^2\Bigg]^{1/2}\label{eq:2006}
\end{align}



\subsection{Numerical Experiment \textbf{$\text{Exp\_EIT}_1$}. Continuous Conductivity.}\label{eit_smooth}

\begin{figure}[h!]
\centering
\includegraphics[scale=0.35, trim=0 0 0 0, clip]{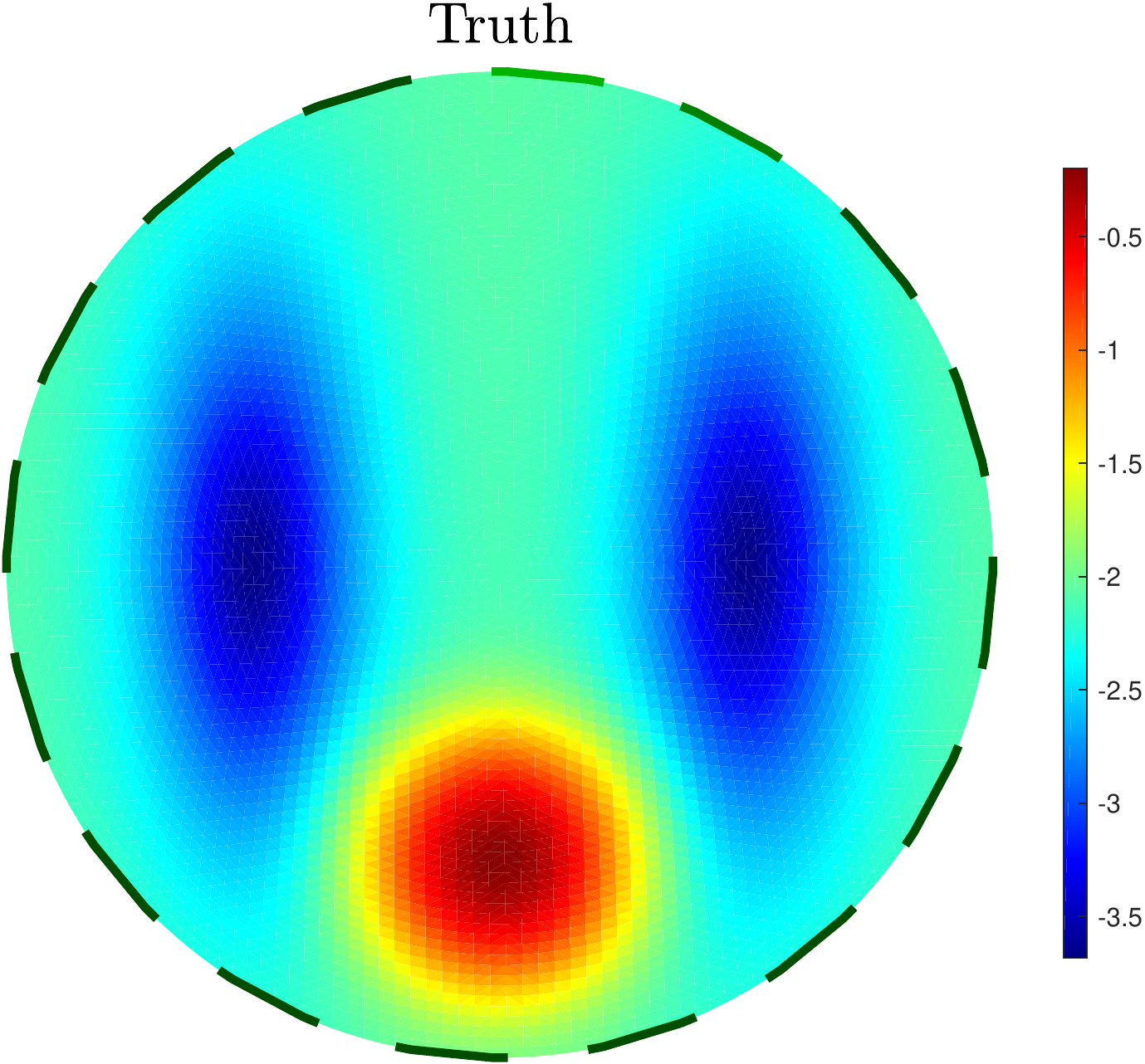} 
\includegraphics[scale=0.35, trim=0 0 0 0, clip]{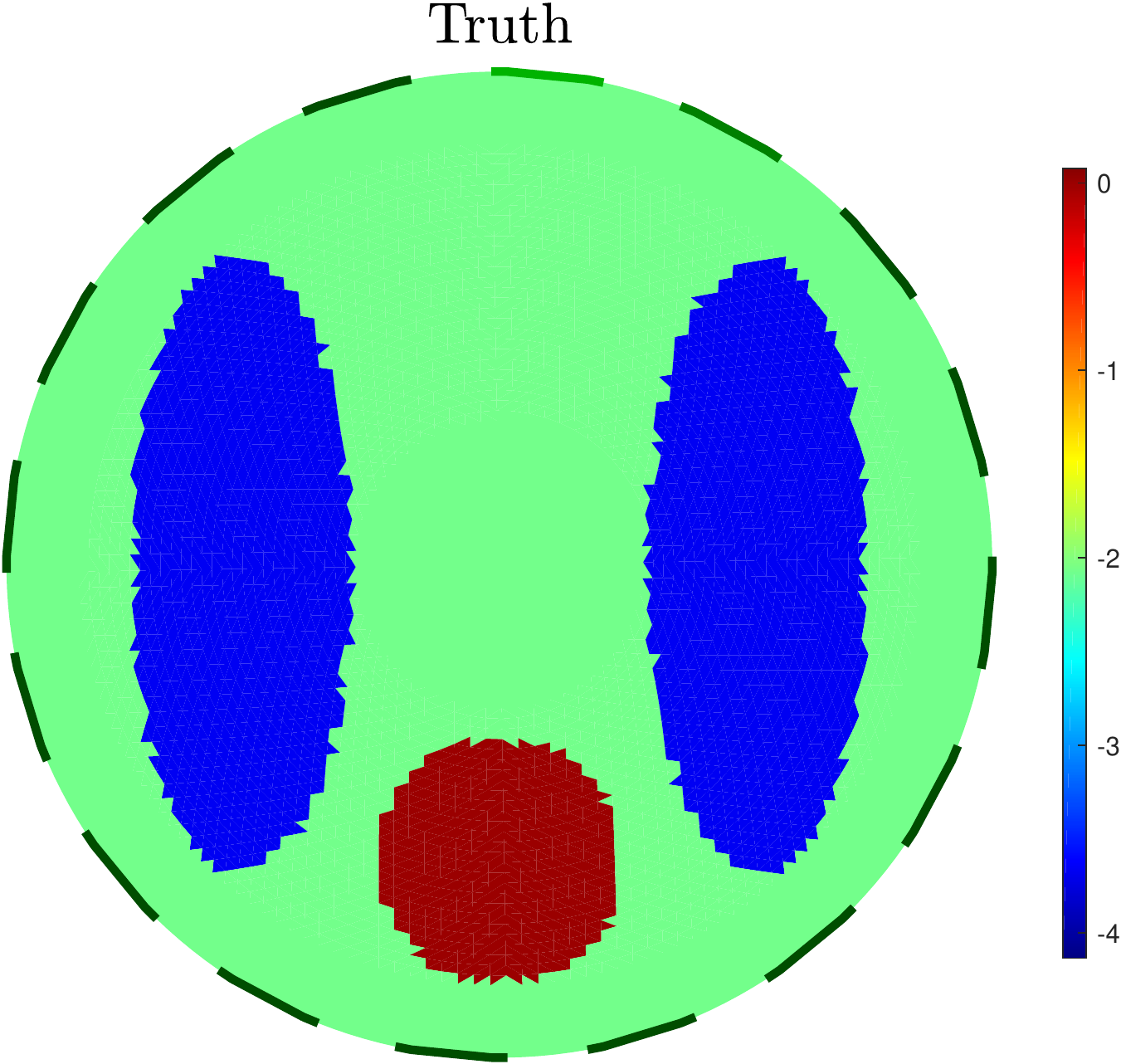}
\includegraphics[scale=0.35, trim=0 0 0 0, clip]{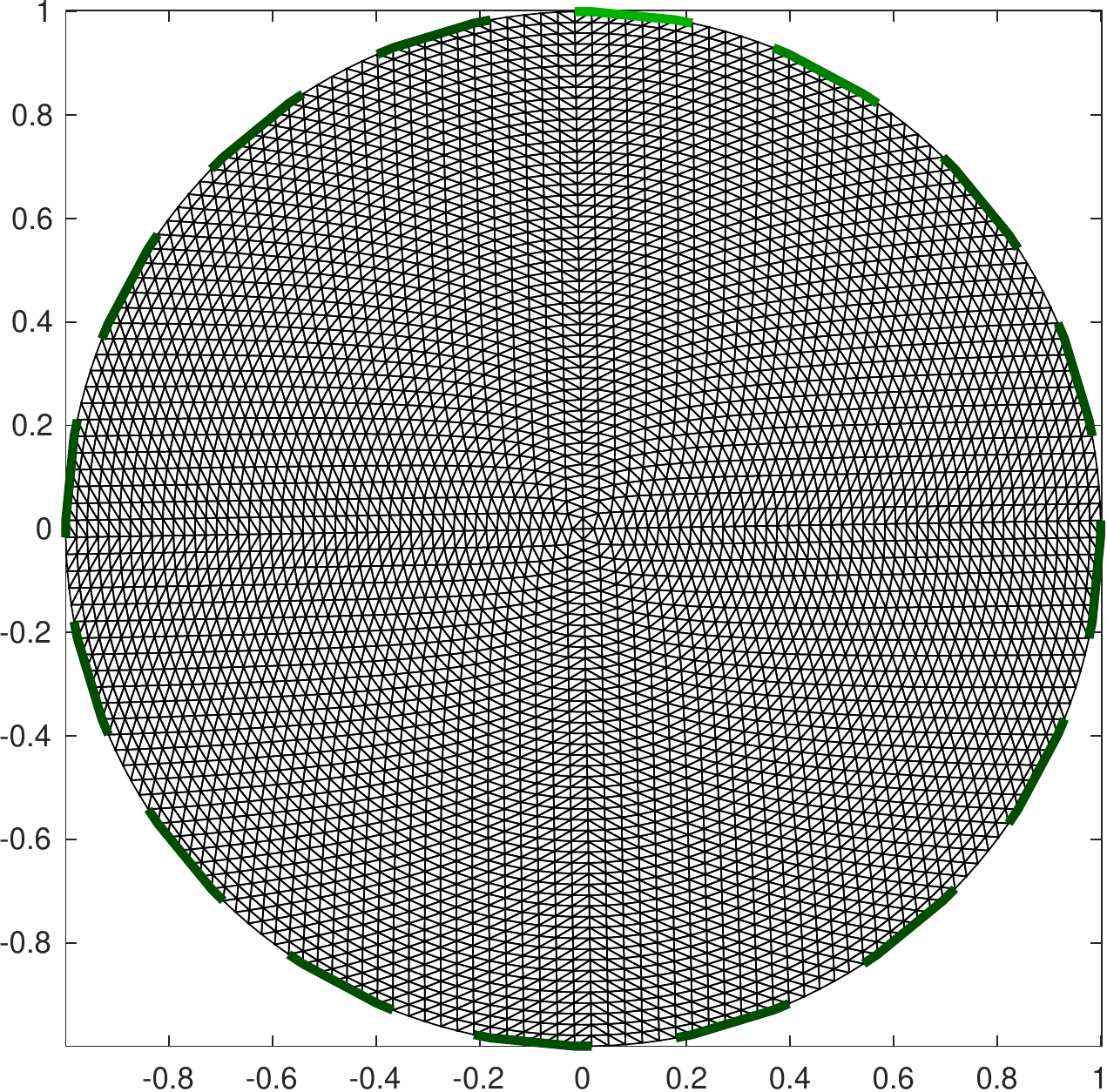}
 \caption{True log conductivity $\log (\kappa^{\dagger})$ for {\bf Exp\_EIT$_1$}  (left) and {\bf Exp\_EIT$_2$} (middle).Right: Mesh for the generation with synthetic data for all EIT experiments.}
\label{Fig1}
\end{figure}

For the first series of experiments the true conductivity, $\kappa^{\dagger}$, is a $C^{\infty}(D)$ function that we specify analytically; the plot of  $\kappa^{\dagger}$ is displayed in Figure \ref{Fig1} (left).  Synthetic data are constructed via $y=V^{\dagger}+\eta$ where $V^{\dagger}=\mathcal{F}(\kappa^{\dagger})$ is computed with the CEM and $\eta$ is a realisation from $N(0,\Gamma)$. We chose $\Gamma=\text{diag}(\gamma_{1},\dots,\gamma_{M})$ where
\begin{align}\label{noi_EIT}
\gamma_{m} = \big(10^{-2}\vert V_{m}^{\dagger}\vert\big)^2 +\big( 10^{-3} \big\vert\max \{ V_{m}^{\dagger}\}_{m=1}^{M}-\min\{V_{m}^{\dagger}\}_{m=1}^{M} \big\vert\Big)^2\qquad m=1,\dots,M.
\end{align}
The first term in the right hand side of (\ref{noi_EIT}) corresponds to adding $1\%$ Gaussian noise. The second term is added simply to avoid small variances from very small voltage (noise-free) measurements.

For this experiment we use the parameterisation of smooth functions, $\kappa=\mathcal{P}_{1}(u)$, from (\ref{para5_2}) with $u=(\lambda, L_{1},L_{2},\omega)$. Note that we have removed $\nu$ and $\sigma$ from the inversion as we discussed in Remark \ref{rem_smooth}; we use fixed values $\nu=3$ and $\sigma=1.5$. The unknown consist of 3 scalars and 1 function, $\omega(x)$, that we discretise on a $100\times 100$ grid. Upon discretisation, the dimension of the unknown is $\text{dim}(\mathcal{U})=10,003$. We follow the discussion of subsection \ref{smoothpara} (see eq. (\ref{para5_5}))  for the selection of the initial ensemble. More specifically, we select $J$ particles $u_{0}^{(j)}=(\lambda^{(j)},L_{1}^{(j)},L_{2}^{(j)},\omega^{(j)})\sim \mu_{0}$ where we define the following prior
\begin{align}\label{rii_EIT}
 \mu_{0}&\equiv U[5\times 10^{-3}, 1]\otimes U[0.15, 0.6]\otimes U[0.15, 0.6]\otimes N(0,\mathbb{I})
 \end{align}
where $U[a,b]$ denotes the uniform distribution on the interval $[a,b]$. 

\begin{figure}[h!]
\centering
\includegraphics[scale=0.46, trim=50 30 0 0, clip]{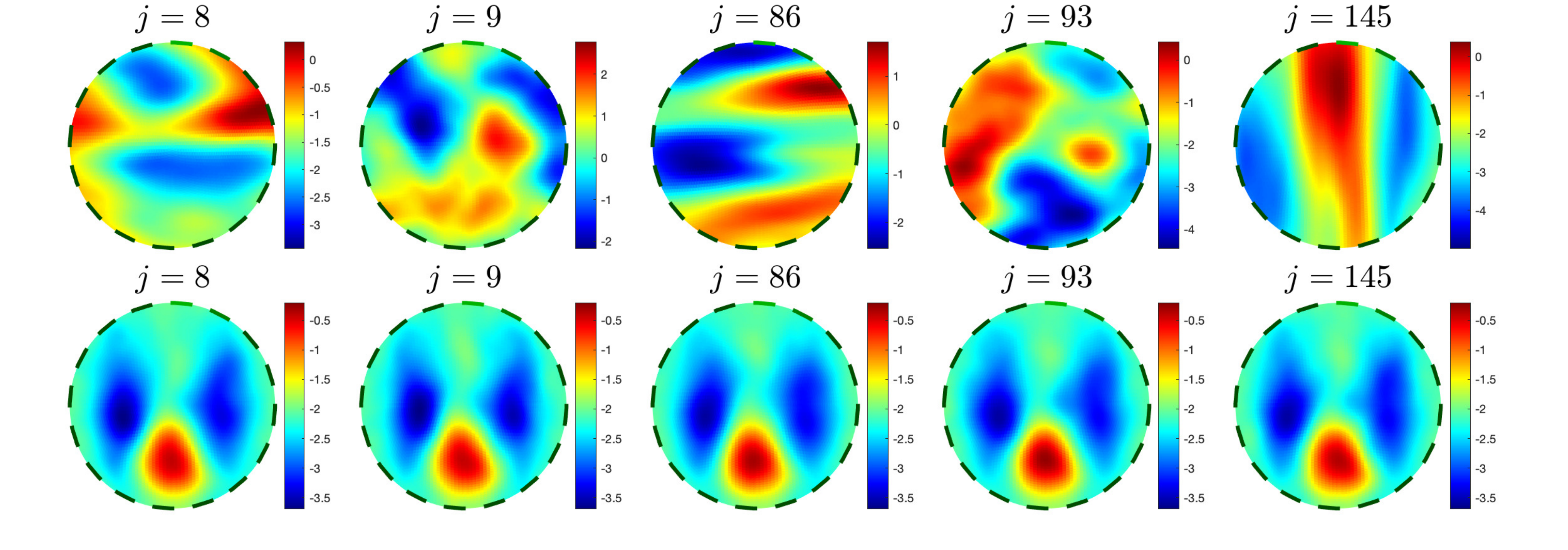}
 \caption{{\bf Exp\_EIT$_1$.}  Top row:  Logarithm of five members from the prior ensemble $\{\kappa_{0}^{(j)}\}_{j=1}^{J}$. Bottom row. Logarithm of five realisations from the final (converged) ensemble  $\{\kappa_{n^*}^{(j)}\}_{j=1}^{J}$.}
\label{Fig2}
\end{figure}

In Figure \ref{Fig2} (top) we show plots for the logarithm of $\kappa_{0}^{(j)}=\mathcal{P}_{1}(u_{0}^{(j)})$ for five of members of the initial ensemble. Note that our choice of smoothness parameter $\nu=3$ produces an initial ensemble that is quite smooth (recall from \ref{smoothpara} that the draws belong to $H^{3-\epsilon}(D)$). 
We also observe substantial differences in the degree of anisotropy which arises from our selection of a reasonably wide distribution of intrinsic lengthscales that we use to produce the initial ensemble.

\begin{figure}[h!]
\centering
\includegraphics[scale=0.335, trim=0 70 0 0, clip]{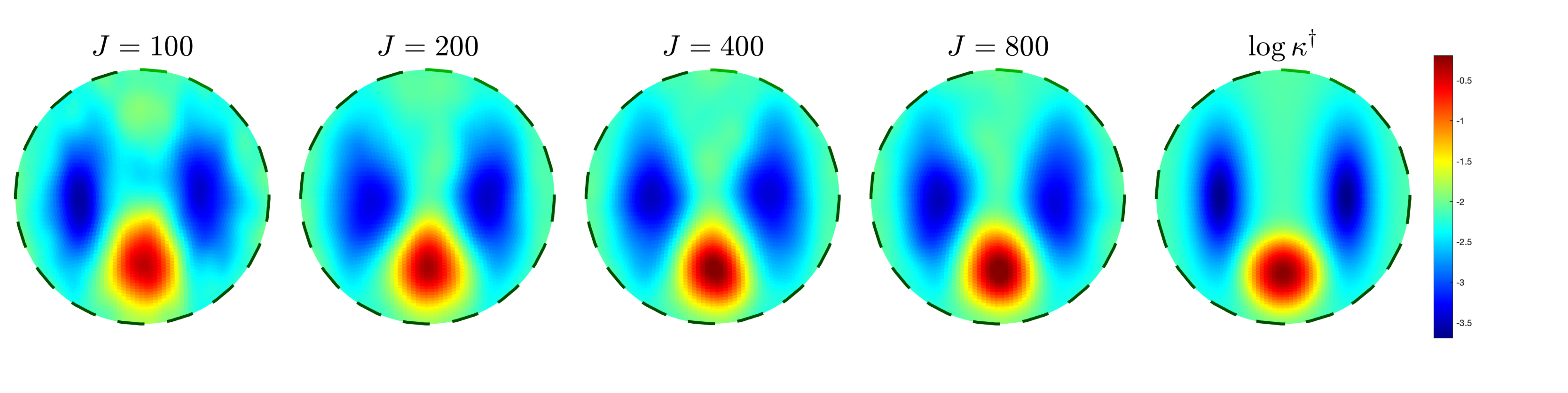}

\caption{{\bf Exp\_EIT$_1$.} Logarithm of $\kappa_{n^*}\equiv\mathcal{P}_{1}(\overline{u}_{n^*})$ computed via the {\bf EKI-DMC} with ensemble size (from left to right) $J=100,\,200,\,400,\,800$. Right panel shows the log of the truth.}
\label{Fig3}
\end{figure}

%

\subsubsection{Results from the inversion using {\bf EKI-DMC} with various choices of $J$.} 
\begin{figure}[h!]
\centering
\includegraphics[scale=0.425, trim=0 0 0 0, clip]{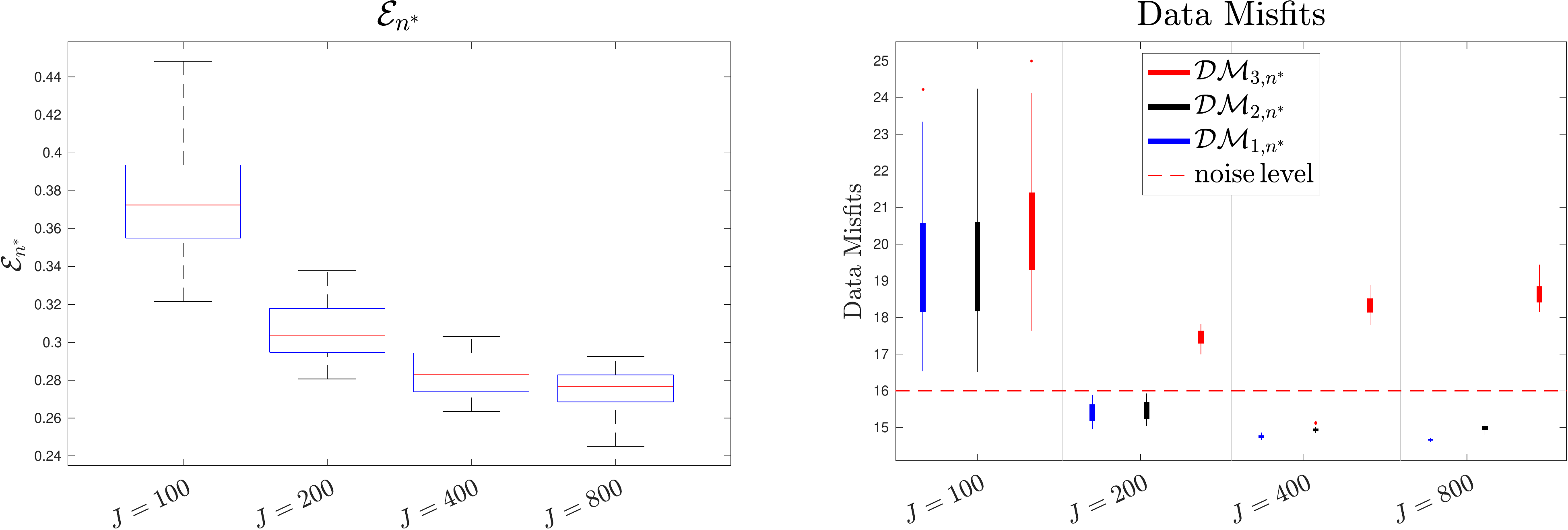} 
\caption{${\bf Exp\_EIT_1.}$ Error with respect to the truth (Left), $\mathcal{E}_{n^{*}}$ (see (\ref{eq:2003})), and data misfits from (\ref{eq:2004})-(\ref{eq:2006}) (Right) computed at the final iteration $n^{*}$ via  {\bf EKI-DMC}. The noise level estimated by $\delta=\sqrt{M}$ is indicated with the dotted red-line in the right panel.}
\label{Fig5}
\end{figure}

Synthetic data and initial ensembles produced as described earlier are used as inputs for {\bf EKI-DMC} (Algorithm \ref{Al2}) with different choices of ensemble size $J$: 100, 200, 400, 800. For each choice of $J$, we conduct 30 experiments with different random selections of the initial ensemble. 
The plots of (log) $\kappa_{n^*}=\mathcal{P}_{1}(\overline{u}_{n^*})$ (i.e. upon convergence) from one of these experiments, computed for each $J$, are displayed in Figure \ref{Fig3} together with the truth (right panel).  In Figure \ref{Fig5} (left) we show boxplots of the relative error w.r.t the truth $\mathcal{E}_{n^*}$ (\ref{eq:2003}), computed at the final iteration, from the set of 30 experiments conducted for each $J$. Boxplots of the data misfits defined in (\ref{eq:2004})-(\ref{eq:2006}) are shown in the right panel of Figure \ref{Fig5} where the mean noise level, approximated via $\delta=\sqrt{M}$, is indicated via the horizontal dotted red line. The average (over the 30 experiments) number of iterations to converge, $n^*$, is displayed in Table \ref{Table1}. 
 These experiments suggest that the choice of $J=200$ provides a reasonable value of the data misfit (i.e. around the noise level). Furthermore, the average relative error w.r.t the truth for $J=200$ (approximately $30\%$) does not improve substantially as we increase $J$. Given these considerations, we select $J=200$ for subsequent runs of  {\bf Exp\_EIT$_1$}. 

 \subsubsection{Further results from one run of {\bf EKI-DMC}  with $J=200$.}\label{eit_smooth_J200} 
\begin{figure}[h!]
\centering
\includegraphics[scale=0.4, trim=50 10 10 10, clip]{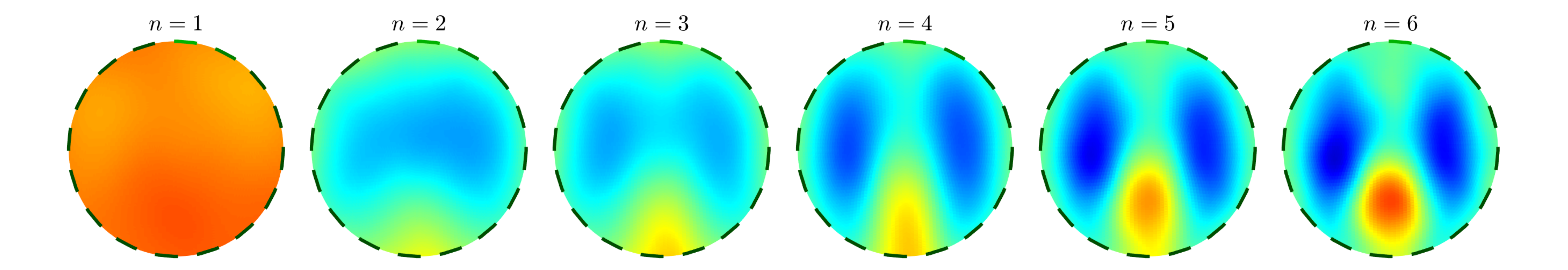} 
\caption{${\bf Exp\_EIT_1.}$ Logarithm of $\kappa_{n}=\mathcal{P}_{1}(\overline{u}_{n})$ computed via the {\bf EKI-DMC} at various intermediate iterations $n$ ($1\leq n\leq n^*$) computed using one ensemble of size $J=200$.}
\label{Fig4}

\end{figure}

\begin{figure}[h!]
\centering
\includegraphics[scale=0.32, trim=0 0 0 0, clip]{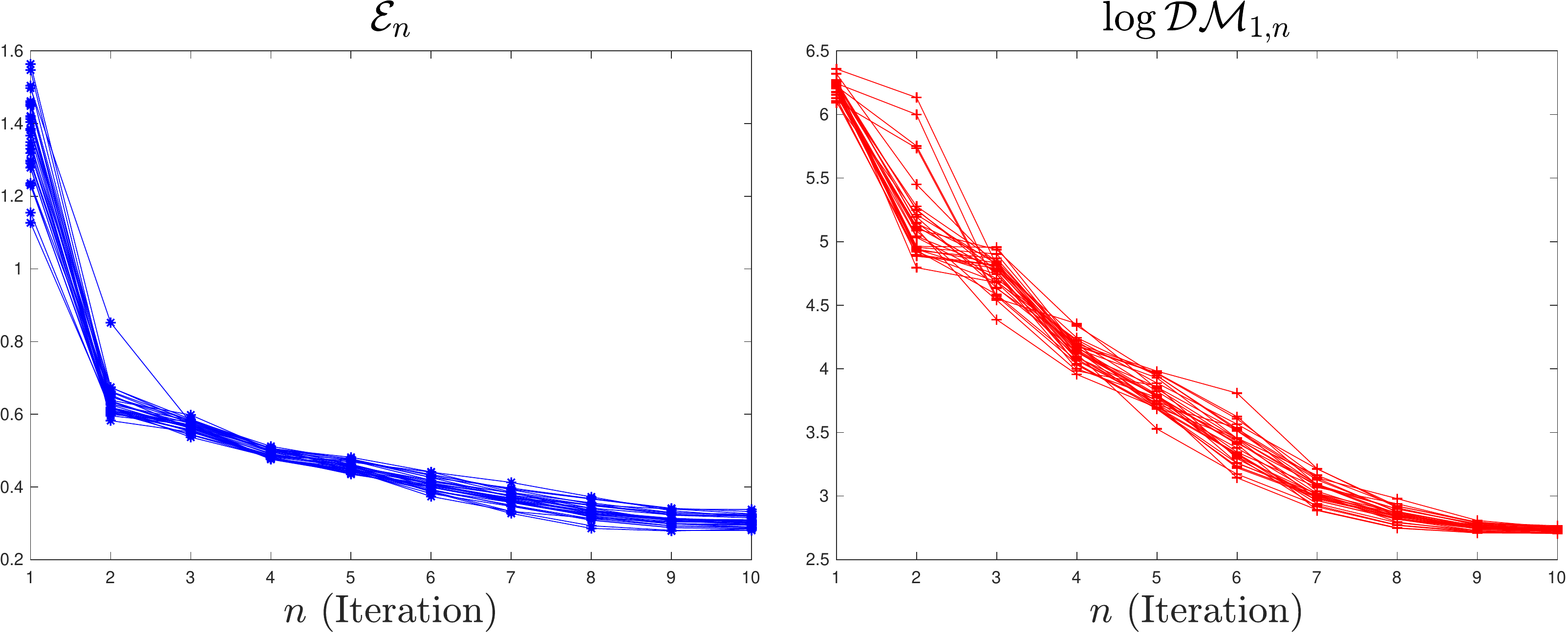} 
\includegraphics[scale=0.3, trim=0 0 0 0, clip]{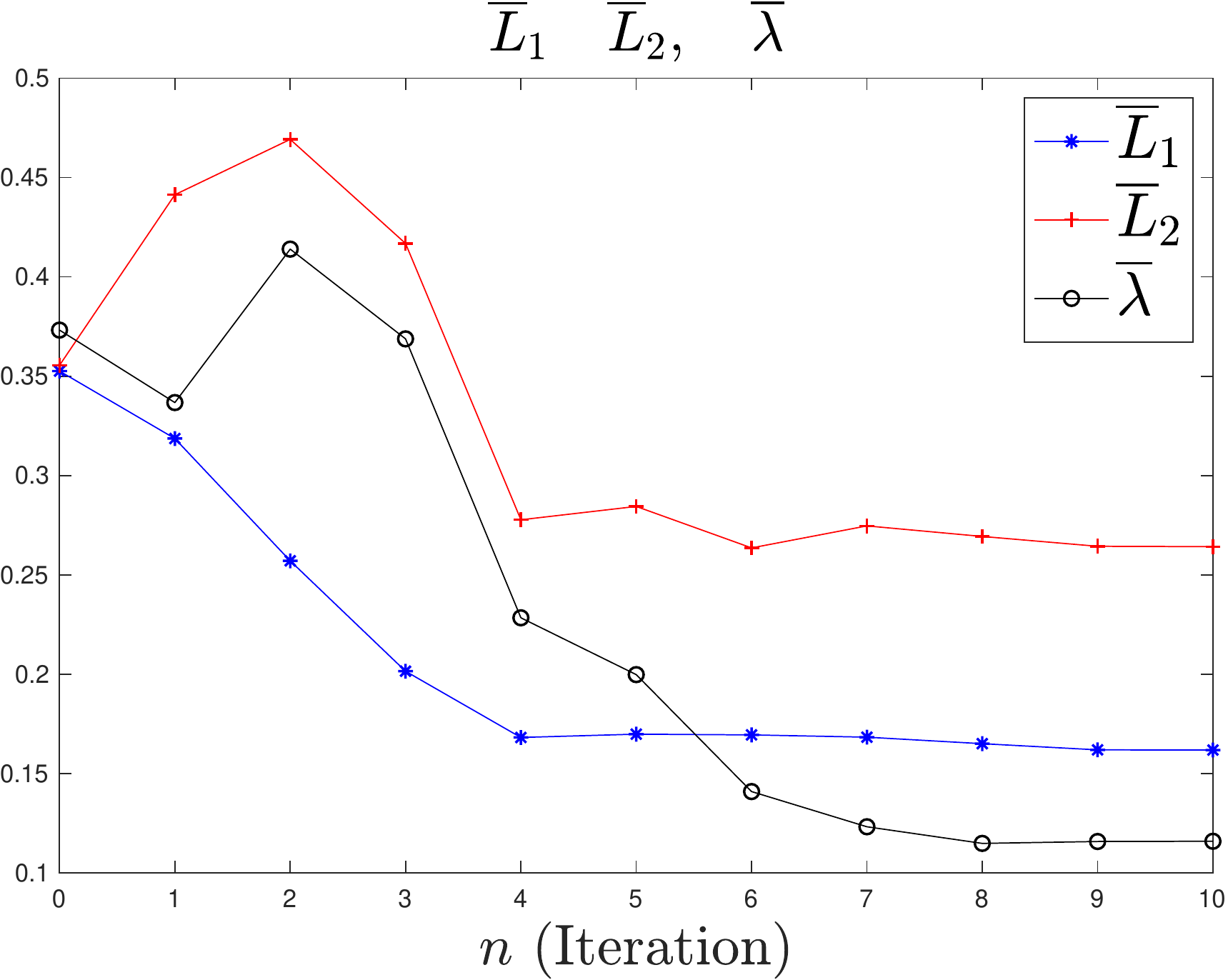}

\caption{${\bf Exp\_EIT_1.}$ Plots of the relative error w.r.t the truth (left), data misfit $\mathcal{DM}_{1,n}$ (middle) and means $\overline{L}_{1}$, $\overline{L}_{2}$ and $\overline{\lambda}$ (right) as a function of $n$. The left and middle panels show results from 30 runs.}
\label{Fig6}
\end{figure}

 In Figure \ref{Fig2} (bottom) we displayed 5 members of the final (converged) ensemble of (log) $\kappa_{n^*}^{(j)}$ corresponding to the initial ensemble from Figure \ref{Fig2} (top). Plots of $\log{\kappa_{n}}$, at some of the intermediate iterations $1\leq n<n^*=10$ can be found in Figure \ref{Fig4}. In Figure \ref{Fig6} (right) we plot, as a function of $n$, the values of the means $\overline{\lambda}$,  $\overline{L}_{1}$, and  $\overline{L}_{2}$, of the ensembles $\{\lambda^{(j)}\}_{j=1}^{J}$, $\{L_{1}^{(j)}\}_{j=1}^{J}$ and $\{L_{2}^{(j)}\}_{j=1}^{J}$, respectively\footnote{We call from (\ref{para5_2})-(\ref{para5_2B}) that $\lambda$, $L_{1}$ and $L_{2}$ are scalar components of the unknown parameter $u$ that we estimate via EKI. For ease in the notation we do not use the subscript $n$ on these variables but we emphasise these are updated at each iteration of EKI}. Note that EKI produces a larger lengthscale in the vertical direction. This comes as no surprise since the truth $\kappa^{\dagger}$ (see Figure \ref{Fig1} (left)) has two inclusions of lower conductivity with larger correlation along the vertical direction.  Finally, in Figure \ref{Fig6} (right and middle) we display, for each of the 30 runs, the relative error w.r.t the truth as well as the (log) data misfit (\ref{eq:2004}), as a function of the iteration number $n$. We note that {\bf EKI-DMC} is very robust and accurate across ensembles. 

\begin{table}[h!]                                                           
\centering                                                                 
\begin{tabular}{|c|c|c|}                                                 
\hline                                                                     
 & {\bf Exp\_EIT$_1$}& {\bf Exp\_EIT$_2$} \\
\hline                                                                     
DM $J=100$ & 10.00$\pm$~0.53& 13.20$\pm$~2.50 \\
\hline                                                                     
DM $J=200$ & 10.00$\pm$~0.00  & 12.83$\pm$~0.59 \\
\hline                                                                     
DM $J=400$ & 10.30$\pm$~0.47  & 14.53$\pm$~0.63 \\
\hline                                                                     
DM $J=800$ & 11.03$\pm$~0.18 & 16.87$\pm$~0.51 \\
\hline       
\end{tabular}                                                              
\caption{Number of iterations $n^{*}$ for EKI using the DMC (\textbf{EKI-DMC}) with various choices of $J$.}   
\label{Table1}  
\end{table}

 \subsubsection{Comparison between  \textbf{EKI-DMC} and  \textbf{EKI-LM}}
 
\begin{figure}[h!]
\centering
\includegraphics[scale=0.42, trim=0 0 0 0, clip]{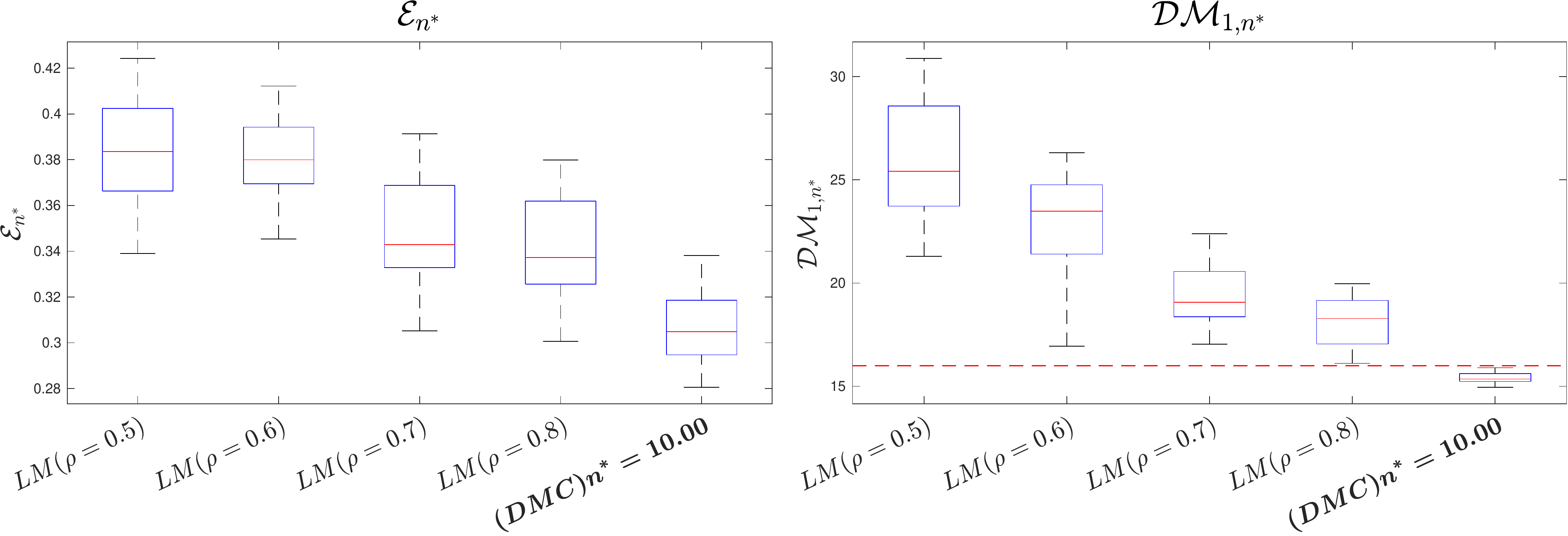} \\

\caption{{\bf Exp\_EIT$_1$.} Error with respect to the truth (left), $\mathcal{E}_{n^{*}}$ (see (\ref{eq:2003})), and data misfit $\mathcal{DM}_{1,n^*}$ (\ref{eq:2004}) (right) computed at the final iteration $n^{*}$ using {\bf EKI-DMC} and {\bf EKI-LM} with various choices of $\rho$. The noise level estimated by $\delta=\sqrt{M}$ is indicated with the dotted red-line in the right panel.}
\label{Fig8}
\end{figure}

\begin{figure}[h!]
\centering
\includegraphics[scale=0.45, trim=60 10 0 7, clip]{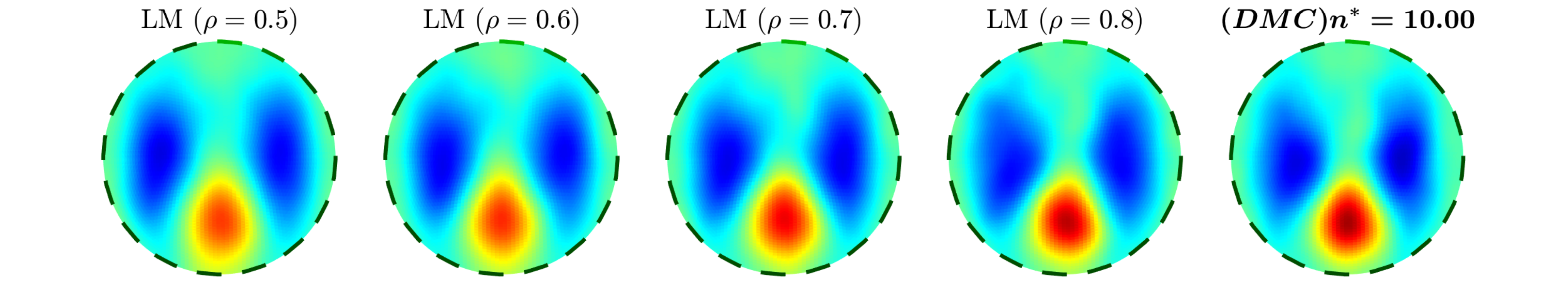}

\caption{{\bf Exp\_EIT$_1$.} Logarithm of $\kappa_{n^*}$ computed using the same initial ensemble ($J=200$) with {\bf EKI-LM} for several choices of  $\rho$. In the top-right panel with display the log of $\kappa_{n^*}$ that we obtain with {\bf EKI-DMC}.}
\label{Fig10}
\end{figure}

We compare the performance of  Algorithms \ref{Al4}-\ref{Al2} using the same set of 30 initial ensembles for each algorithm using $J=200$. In \cite{Neil,Iglesias2016}, this selection of $J$ was sufficient to provide stable and accurate estimates for EIT. We consider different choices of the input $\rho$ in \textbf{EKI-LM} and, for simplicity, we set $\tau=1/\rho+10^{-6}$. In Figure \ref{Fig8} we show boxplots of the error w.r.t the truth (right) and the data misfit $\mathcal{DM}_{1,n^*}$ (left) obtained with several choices of $\rho$ for \textbf{EKI-LM}; the results from \textbf{EKI-DMC} for $J=200$ are also included in these plots.  Similar behaviour is observed for $\mathcal{DM}_{2,n^*}$ and  $\mathcal{DM}_{3,n^*}$ and so these plots are omitted. The number of iterations for \textbf{EKI-LM} to achieve convergence, $n^*$, is displayed in Table \ref{Table1B}. For one of the 30 runs, in Figure \ref{Fig10} we show the plots of $\log{\kappa_{n^*}}$ computed with both algorithms using the selections of parameters described above.

While \textbf{EKI-DMC} does not depend on any tuning parameter, the results above show that \textbf{EKI-LM} is highly dependent on the selection of $\rho$ (which is specified a priori). From Table \ref{Table1B} we see that, for $\rho=0.6$, the computational cost of \textbf{EKI-LM} is similar to the cost of \textbf{EKI-DMC}; i.e. convergence in approximately 10 iterations (in average). For this $\rho$, the average error computed with \textbf{EKI-LM} ($\approx 38\%$) is, however, larger than the error obtained via \textbf{EKI-DMC} ($\approx 0.305\%$). Although the improvement in accuracy of \textbf{EKI-DMC} over \textbf{EKI-LM} may not be overly impressive, from Figure \ref{Fig10} (computed from one run) we note that the area of high conductivity is not accurately captured by  \textbf{EKI-LM} with $\rho<0.8$.  For these experiments, if instead of $\rho=0.6$ we choose, say $\rho=0.8$ (see Table \ref{Table1B}), the computational cost of \textbf{EKI-LM} doubles without improving its accuracy with respect to \textbf{EKI-DMC}.

Furthermore, from Figure \ref{Fig8} we also observe that both the data-misfit $\mathcal{DM}_{1,n^*}$ and the error w.r.t the truth obtained via \textbf{EKI-LM} decreases as we increase $\rho$. This behaviour is expected since large $\rho$ implies smaller $\tau$ and, hence, smaller data misfits upon convergence (see stopping criteria for \textbf{EKI-LM} in eq. (\ref{EQ2})). Even though the measures of performance of \textbf{EKI-LM} (see Figure \ref{Fig8}) seem to approach those of \textbf{EKI-DMC} for increasing $\rho$, experiments (not shown) conducted with even larger $\rho$'s, (i.e. $\rho\ge 0.9$) do not yield the same level of accuracy that we achieve with \textbf{EKI-LM}. This is due to the fact that the stopping criteria in \textbf{EKI-LM} does not allow the algorithm to iterate once $\mathcal{DM}_{1,n^*}$ is smaller than $\tau\delta=\tau \sqrt{M}$ (recall $\tau>1/\rho$). Hover, for this particular problem, it seems that iterating below this threshold yields more accurate estimates. Indeed, we note that for the estimates obtained by \textbf{EKI-DMC}, the value of $\mathcal{DM}_{1,n^*}$ is smaller than the $\tau\delta$ for all the 30 runs. We conclude that, for this particular problem setting, the stopping criteria in \textbf{EKI-LM} is not optimal. 

\begin{figure}[h!]
\centering
\includegraphics[scale=0.55, trim=0 0 0 0, clip]{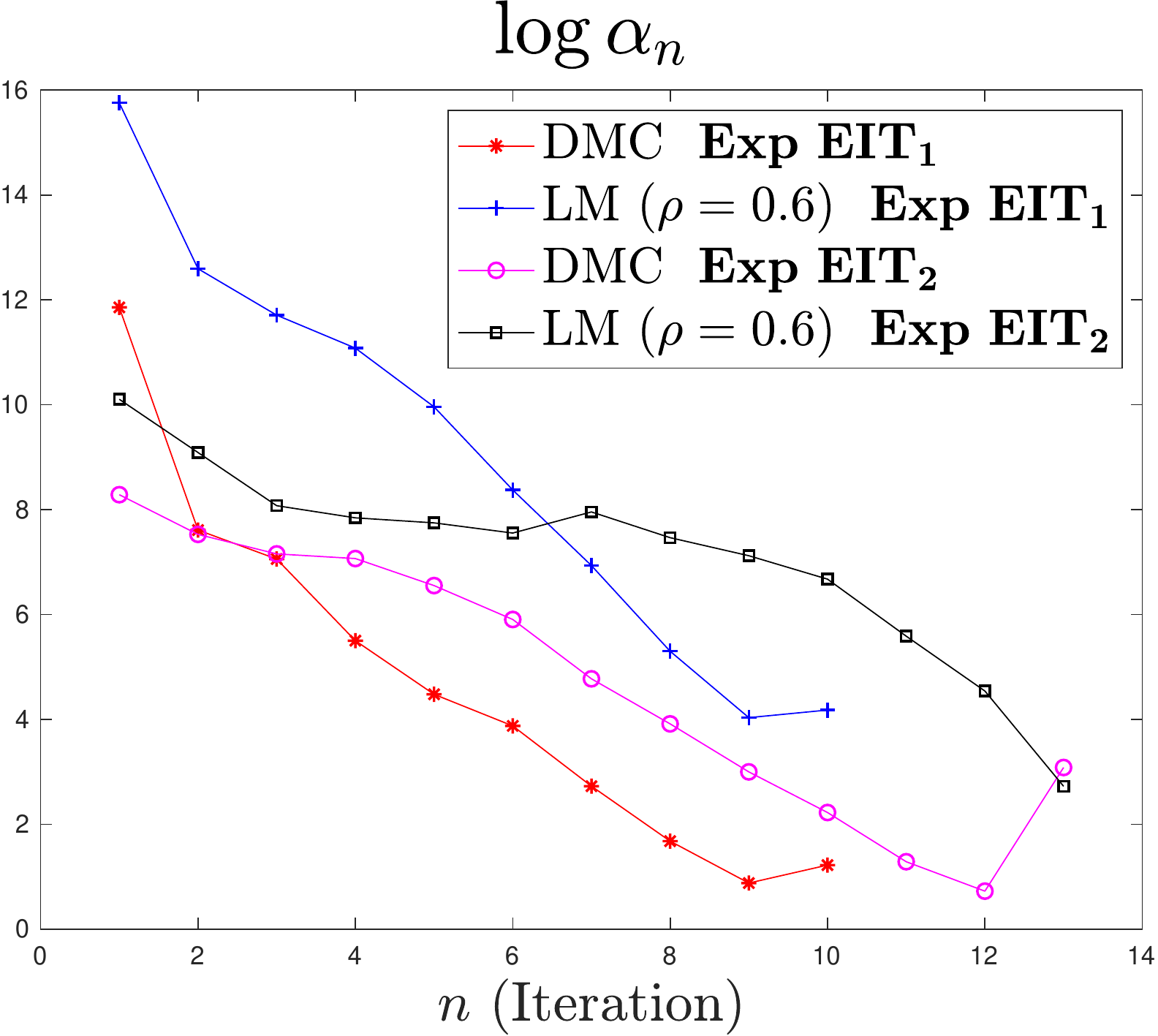}

\caption{Logarithm of the regularisation parameter $\alpha_{n}$ as a function of $n$ computed with {\bf EKI-DMC} and {\bf EKI-LM} for EIT experiments. }
\label{Fig9}
\end{figure}

In order to further understand the potential limitation of the stopping criteria in \textbf{EKI-LM} as proposed in \cite{EnsembleYo}, in Figure \ref{Fig9} we display the behaviour of $\log \alpha_{n}$ as a function of $n$, computed from one run (with same initial ensemble) of \textbf{EKI-DMC} and \textbf{EKI-LM} with $\rho=0.6$. As mentioned above, for this particular selection of $\rho$, both \textbf{EKI-LM}  and \textbf{EKI-DMC} converged in 10 iterations so the cost of running both algorithms is the same. We find that the selection of $\alpha_{n}$ via \textbf{EKI-LM} yields quite large values of $\alpha_{n}$ compared to those obtained via \textbf{EKI-DMC}. Larger $\alpha_{n}$'s, as discussed in subsection \ref{new_lab}, means slower convergence. Although we can obtain smaller $\alpha_{n}$'s using \textbf{EKI-LM} with smaller $\rho$'s, this will result in larger $\tau$'s and, consequently, earlier termination with is detrimental to the accuracy of the algorithm. 

The aforementioned limitation of the stopping criteria in \textbf{EKI-LM} can be addressed, for example, by adopting the same stopping rule that we use in \textbf{EKI-DMC}. In other words, we may replace (\ref{EQ2}) by imposing that the $\alpha_{n}$'s in \textbf{EKI-LM} satisfy (\ref{eqA:9}). Such an approach has been used, for example, in \cite{other} in the context of history matching of petroleum reservoirs. We have repeated our experiments (not shown) using (\ref{eqA:9}) as stopping rule for \textbf{EKI-LM} which show that, indeed, using (\ref{eqA:9}) enable us to run the algorithm beyond the noise level and, hence, obtain more accurate estimates at a lower computational cost. However, the role of $\rho$ still determines the performance of \textbf{EKI-LM} and, while (\ref{eqA:9}) allows us to use smaller $\rho's$ (thus smaller $\alpha_{n}$'s and faster convergence), for some of these small values the algorithm loses stability.

 In summary, we may find problem-specific tuning parameters $\rho$ and $\tau$ in \textbf{EKI-LM}, or modifications of \textbf{EKI-LM} by adopting (\ref{eqA:9}) as stopping rule, that will yield comparable level of accuracy and even similar computational cost to those achieved by  \textbf{EKI-DMC}. However, finding optimal stopping criteria and/or tuning parameters may require thorough numerical testing which is often computationally intensive and unsuitable for practical large-scale applications. 
 

\begin{table}[h!]                                                           
\centering                                                                 
\begin{tabular}{|c|c|c|}                                                 
\hline                                                                     
 &{\bf Exp\_EIT$_1$}& {\bf Exp\_EIT$_2$}\\
\hline                                                                     
LM $\rho=0.5$ & 8.57$\pm$~0.50  & 10.30$\pm$~0.79 \\       
\hline                                                                     
LM $\rho=0.6$ & 10.53$\pm$~0.51 & 14.10$\pm$~1.79 \\      
\hline                                                                     
LM $\rho=0.7$ & 14.40$\pm$~0.56  & 20.03$\pm$~1.45 \\     
\hline                                                                     
LM $\rho=0.8$ & 21.80$\pm$~0.76 & 32.53$\pm$~2.65 \\     
\hline       
\end{tabular}                                                              
\caption{Number of iterations $n^{*}$ for EKI using the LM approach (\textbf{EKI-LM}) with various choices of $\rho$.}   
\label{Table1B}  
\end{table}

\subsection{Numerical Experiment \textbf{$\text{Exp\_EIT}_2$}. Discontinuous (piecewise constant) conductivity.}\label{eit_disc}

The true conductivity for the experiments of this section is the piecewise constant function with plot displayed in Figure \ref{Fig1} (middle). The background, low and high conductivity regions have constant values:
\begin{align}\label{para5_12B00}
\kappa_{b}^{\dagger}=0.125,\qquad \kappa_{l}^{\dagger}=0.025,\qquad \kappa_{h}^{\dagger}=1.0,
\end{align}
respectively.  Synthetic data are generated in the same way as described in the previous subsection, and with noise covariance from (\ref{noi_EIT}).

\begin{figure}[h]
\centering
\includegraphics[scale=0.3, trim=0 70 0 25, clip]{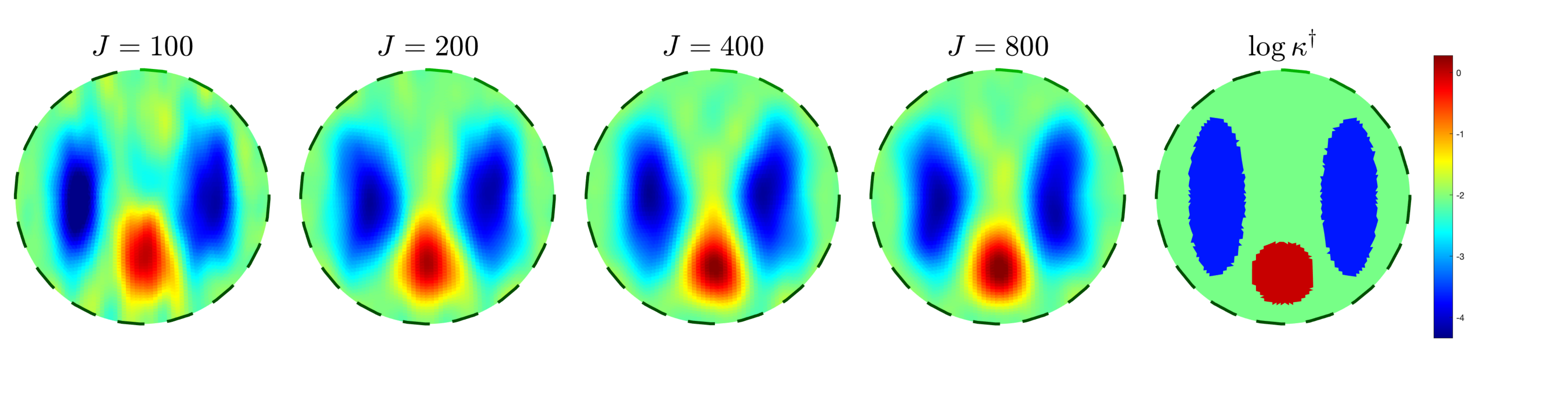} 
\caption{Logarithm of $\kappa_{n^*}$ computed via {\bf EKI-DMC} ensemble size (from left to right) $J=100,\,200,\,400,\,800$. Right panel shows $\log (\kappa^{\dagger})$. This results use the parameterisation from ${\bf \text{Exp\_EIT}_1}$.}
\label{Fig_ape1}
\end{figure}

In Figure \ref{Fig_ape1} we show the plots of (log) $\kappa_{n^*}$ computed from one run with {\bf EKI-DMC} using the same parameterisation (for continuous functions) and initial ensemble from \textbf{$\text{Exp\_EIT}_1 $}. While these results show that we can identify the three different regions of different conductivity, the reconstruction of the interface between these regions is quite inaccurate because of smoothness enforced by the parameterisation used within EKI. We overcome this limitation by means of  a level-set parameterisations within the EKI framework. 

\subsubsection{Level-set parameterisation and prior ensemble}

For the series of experiment in this subsection we test EKI with the level-set parameterisation $\kappa=\mathcal{P}_{2}(u)$ from (\ref{para5_13}). As discuss in Remark \ref{rem_smooth} we remove $\sigma_{f}$ and $\nu_{f}$ from the inversion. Here we use the values $\nu_{f}=2$ and $\sigma_{f}=0.5$. In addition, we also leave $\lambda_{f}$ out of the inversion; we take $\lambda_{f}=1$. This selection is made for simplicity since we expect to capture all the variability of the level-set function via the term $\mathcal{W}_{\Theta_{f}}\omega_{f}$ in  (\ref{para5_11B}). The unknown reduces to 
$$u= (\kappa_{l},\kappa_{b},\kappa_{h},L_{1,f},L_{2,f},\omega_{f}),$$
which  consists of 5 scalars and one function, $\omega_{f}$, which we discretise also on a $100\times 100$ grid. Hence, the dimension of the unknown $u$ is $\text{dim}(\mathcal{U})=10,005$. We select the initial ensemble according to
\begin{eqnarray}\label{para5_12BBB}
u_{0}^{(j)}\sim U[0.015, 0.075]\otimes U[0.1, 0.4] \otimes U[0.65, 1.1]\otimes U[0.15, 0.6]\otimes U[0.15, 0.6]\otimes N(0,\mathbb{I}),
\end{eqnarray}
From (\ref{para5_11B}) we know that this selection produces a centred ensemble of initial level-set functions. It is worth noticing from (\ref{para5_12BBB}) that there is no overlapping among the support of the distributions for the conductivity values (i.e. $\kappa_{l}$, $\kappa_{b}$ and $\kappa_{h}$) on each region. Furthermore, each of these intervals contain the true values form (\ref{para5_12B00}). Therefore, we work under the assumption that (i) there is clear contrast between the unknown values on each region and (ii) we have good knowledge of the range of possible values for the unknown conductivity in each of those regions.

\begin{figure}[h!]
\centering
\includegraphics[scale=0.3, trim=0 20 0 0, clip]{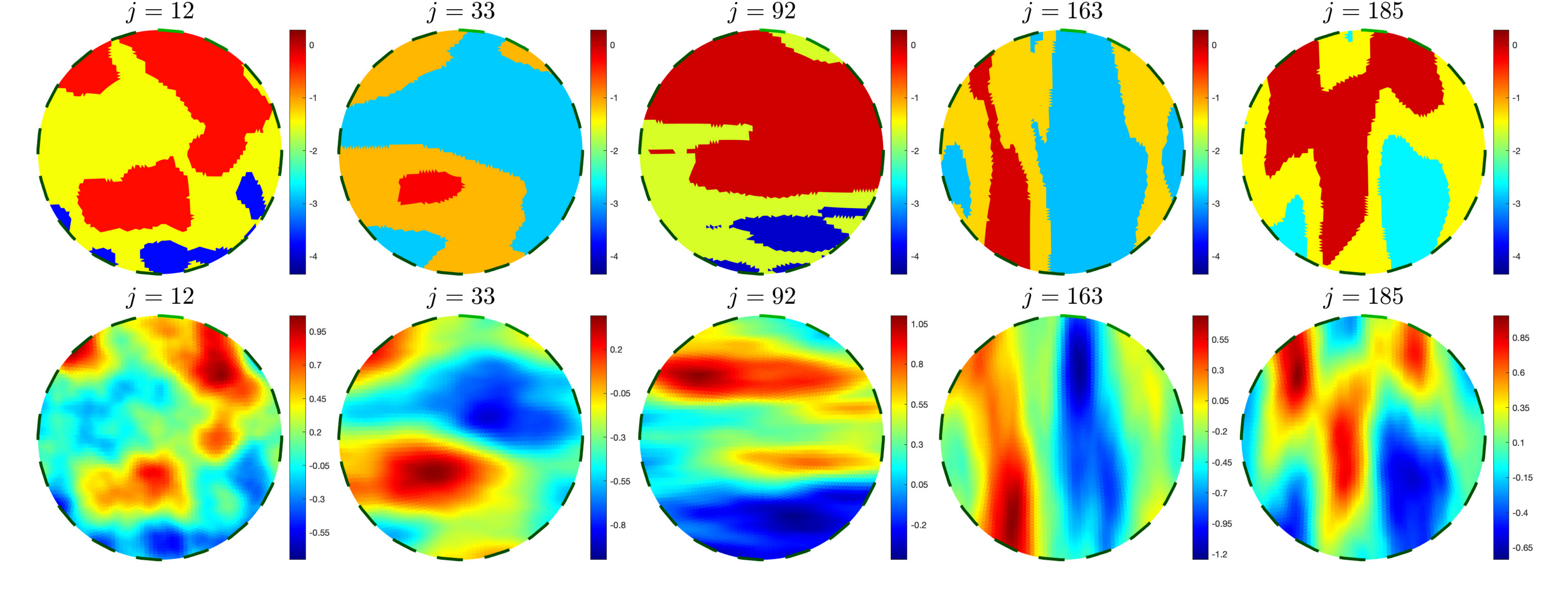} 

\caption{${\bf Exp\_EIT_2}$. Five members of the Initial  ensemble  of (log) $\{\kappa_{0}^{(j)}\}_{j=1}^{J}$ (top row) and their corresponding level-set set function $\{f_{0}^{(j)}\}_{j=1}^{J}$ (bottom row).}
\label{Fig17}
\end{figure}

In Figure \ref{Fig17} we show the plots of some $\log \kappa_{0}^{(j)}=\log{\mathcal{P}_{2}(u_{0}^{(j)})}$ (top panels) and the corresponding level-set functions $f_{0}^{(j)}=\log (\mathcal{P}_{1}(u_{f,0}^{(j)}))$ (bottom panels) from five ensemble members. Our selection $\nu_{f}=2$ imposes smoothness in the level-set function and, in turn, in the interface between the three regions of low, high and background conductivity. From Figure \ref{Fig17} (bottom) we notice that the ensemble of level-sets displays anisotropy induced by randomising the intrinsic lengthscales in the level-set function. This variability can be seen in the corresponding interface between regions of different conductivities (top row Figure \ref{Fig17}). Note that the values of conductivity within each region are variable across particles.

\subsubsection{Results from several choices of $J$ in  \textbf{EKI-DMC}.} 

\begin{figure}[h!]
\centering
\includegraphics[scale=0.42]{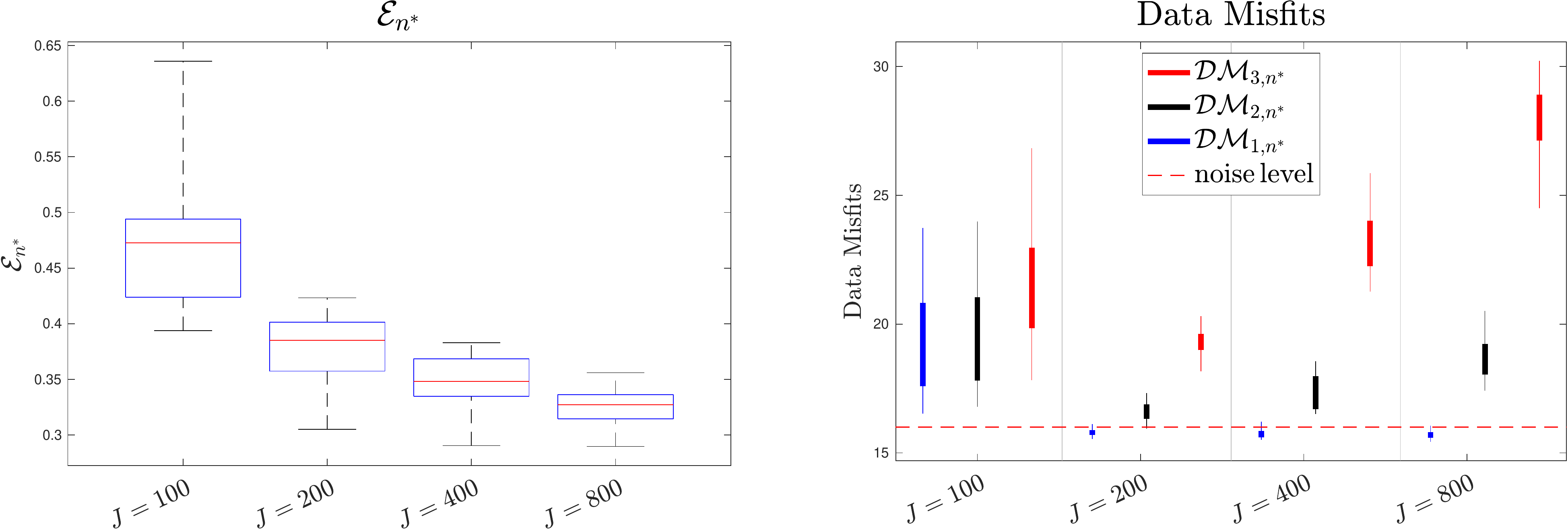}

\caption{${\bf Exp\_EIT_2}$. Error with respect to the truth (Left), $\mathcal{E}_{n^{*}}$ (see (\ref{eq:2003})), and data misfits from (\ref{eq:2004})-(\ref{eq:2006}) (Right) computed at the final iteration $n^{*}$ with {\bf EKI-DMC} with different choices of $J$. The noise level estimated by $\delta=\sqrt{M}$ is indicated with the dotted red-line in the right panel.}
\label{Fig11}
\end{figure}

\begin{figure}[h!]
\centering
\includegraphics[scale=0.335, trim=10 50 0 20, clip]{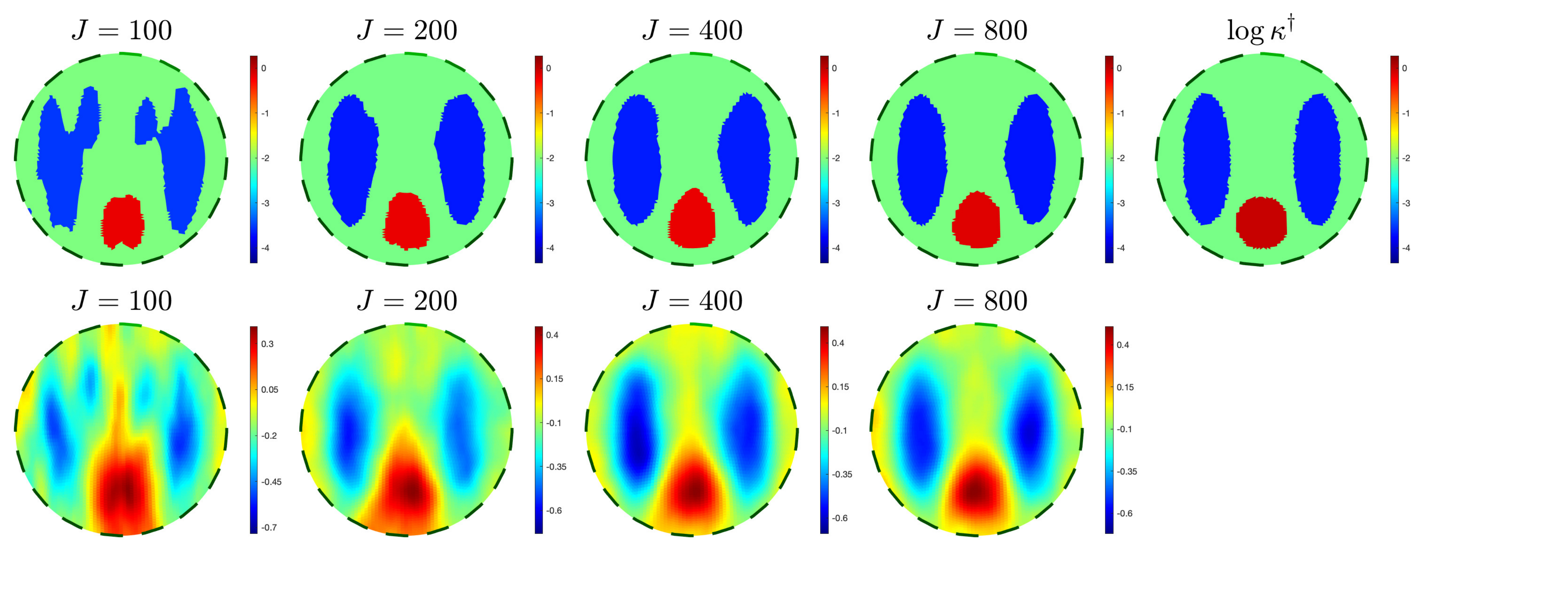}

\caption{${\bf Exp\_EIT_2.}$ Plots of (log) $\kappa_{n^*}$ (top row) and the corresponding level-set function $f_{n^*}$ (bottom row) computed with {\bf EKI-DMC} with different choice of  ensemble size $J$. Top-right panel shows $\log (\kappa^{\dagger})$.}

\label{Fig14}
\end{figure}

In Figure \ref{Fig11} we show boxplots of the relative error w.r.t the truth as well as data misfits (\ref{eq:2004})-(\ref{eq:2006}). As before, these are results from 30  \textbf{EKI-DMC} runs using different initial ensembles for each choice of $J$. We can clearly see a decrease in the error w.r.t the truth as we increase $J$, while the data misfits $\mathcal{DM}_{1,n}$ and $\mathcal{DM}_{2,n}$ achieve values close to the noise level (indicated via dotted red line) for $J\ge 200$. We see that, again, $J=200$ is a good compromise between accuracy and cost. Using $J=400$ will double the cost with a marginal improvement in accuracy. The choice of $J=200$ also yields reasonable visual results as we can verify from Figure \ref{Fig14} where, for one run of \textbf{EKI-DMC}, we display the log of $\kappa_{n^{*}}=\mathcal{P}_{2}(\overline{u}_{n^*})$ (top panels) and the level-set function $f_{n^*}=\log (\mathcal{P}_{1}(\overline{u}_{f,n^*}))$ (bottom panels).


Although we note that the average data misfit $\mathcal{DM}_{3,n^*}$ in Figure \ref{Fig11} seems to increase with $J$, experiments (not shown) with larger $J$ suggest that $\mathcal{DM}_{3,n^*}$ eventually stabilises. The increase in $\mathcal{DM}_{3,n^*}$ can be attributed to the fact that larger $J$'s produces a better spread/coverage of the distribution of particles; some of these particles yield a larger data misfit within EKI\footnote{from the Bayesian perspective some of these particles have small (approximate) posterior probability}. The ensemble mean $\overline{u}_{n}$, however, is quite accurate (see Figure. \ref{Fig14}, for $J=200$) and so the corresponding $\kappa_{n}=\mathcal{P}_{2}(\overline{u}_{n})$ yields reasonable values of the $\mathcal{DM}_{2,n^* }$ (i.e. close to the noise level). 

 \subsubsection{Results from one run of \textbf{EKI-DMC}  with $J=200$.}
\begin{figure}[h!]
\centering
\includegraphics[scale=0.3,, trim=0 30 0 0, clip]{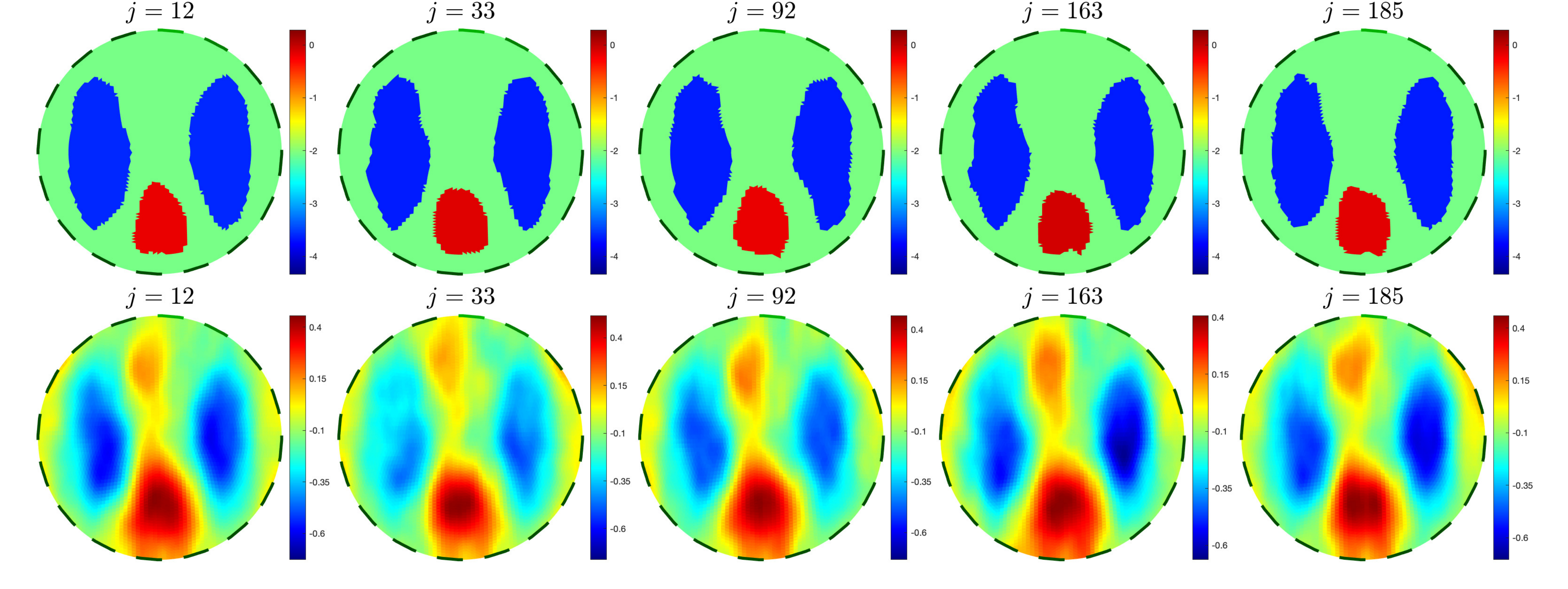}

\caption{${\bf Exp\_EIT_2.}$ Five members of the converged ensemble of (log) $\{\kappa_{n^*}^{(j)}\}_{j=1}^{J}$ (top row) and their corresponding level-set set function $\{f_{n^*}^{(j)}\}_{j=1}^{J}$ (bottom row).}
\label{Fig18}
\end{figure}

\begin{figure}[h!]
\centering
\includegraphics[scale=0.42]{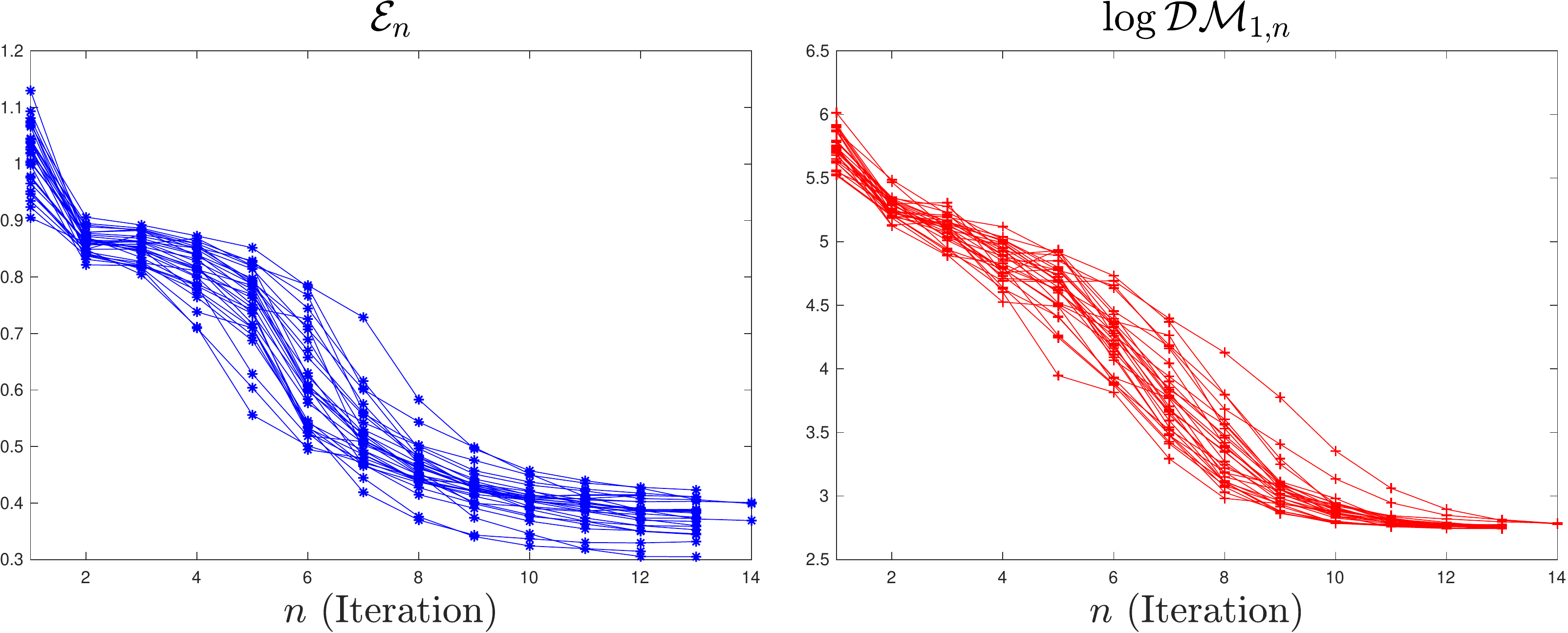} 
\caption{${\bf Exp\_EIT_2.}$ Plots of the relative error w.r.t the truth (left) and data misfit $\mathcal{DM}_{1,n}$ (right) computed from 30 runs with {\bf EKI-DMC}.}

\label{Fig15C}
\end{figure}

Figure \ref{Fig18} shows some members from the final (converged) ensemble corresponding to the initial ensemble from Figure \ref{Fig17}; the top panels are $\log{\kappa_{n^*}^{(j)}}=\mathcal{P}_{2}(u_{n^*}^{(j)})$ while the bottom panels show the level-set functions $f_{n^*}^{(j)}=\log (\mathcal{P}_{1}(u_{f,n^*}^{(j)}))$. This figure shows that there is, indeed, significant variability across particles of the converged ensemble $\kappa_{n^*}^{(j)}$ and, in turn, possible large spread in the values of the data misfit obtained using each of these conductivities (hence potentially large values of $\mathcal{DM}_{3,n}$). The average error w.r.t. the truth and (log) data misfit $\mathcal{DM}_{1,n}$ are shown in Figure \ref{Fig15C}, as a function of the iteration number $n$ (these are results from our 30 runs). Note that, in contrast to the previous experiment, the error displays more variability across ensembles.

\begin{figure}[h!]
\centering
\includegraphics[scale=0.35, trim=10 20 10 10, clip]{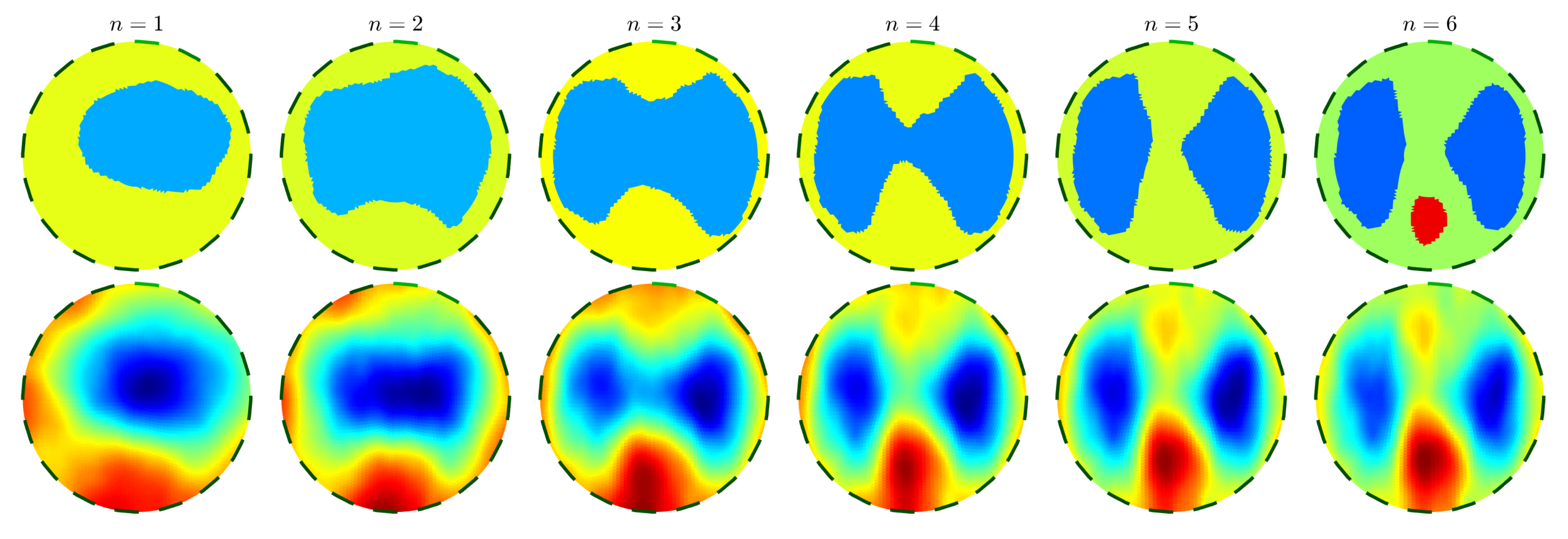} 
\caption{${\bf Exp\_EIT_2.}$ Logarithm of $\kappa_{n})$ (top) and the corresponding level-set $f_{n}$ (bottom) computed with one run of {\bf EKI-DMC} (with $J=200$) at various intermediate iterations $n$ ($1\leq n\leq n^*$).}
\label{Fig15}
\end{figure}

\begin{figure}[h!]
\centering
\includegraphics[scale=0.33, trim=0 0 0 0, clip]{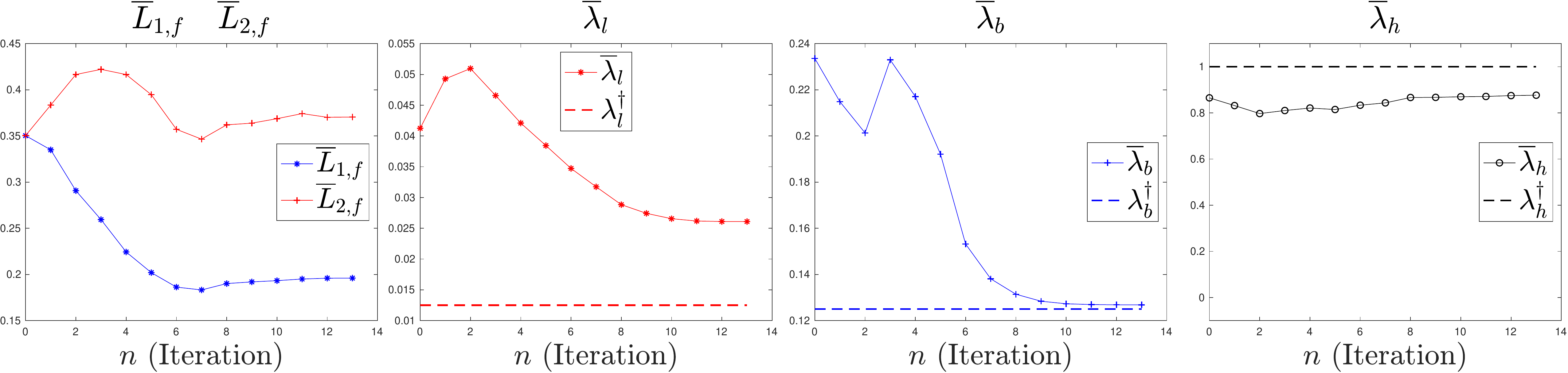} 
\vskip2mm
\includegraphics[scale=0.42, trim=0 0 0 0, clip]{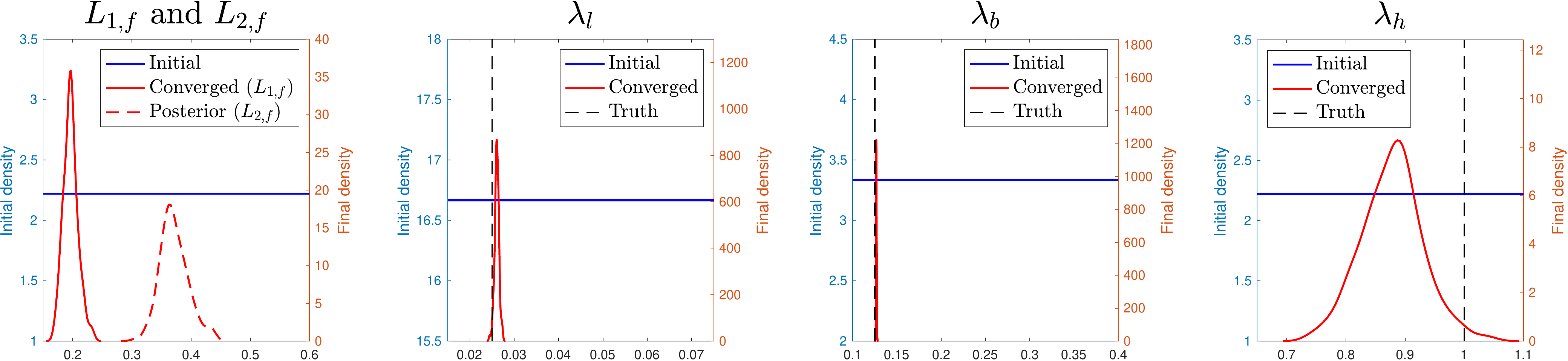}
\caption{${\bf Exp\_EIT_2.}$ Top: Plots of the means $\overline{L}_{1,f}$ and $\overline{L}_{2,f}$ (left) as well as the mean conductivity values $\overline{\lambda}_{l}$ (left-middle), $\overline{\lambda}_{b}$ (middle-right) and $\overline{\lambda}_{h}$ (right) corresponding to the low, background and high conductivity regions. Bottom: densities from the initial and final (converged) ensemble of $L_{1,f}$, $L_{2,f}$, $\lambda_{l}$, $\lambda_{b}$, $\lambda_{h}$.}

\label{Fig15B}
\end{figure}

Plots of $\log{\kappa_{n}}$ and the level-set function $f_{n}$, at some of the intermediate iterations $1\leq n<n^*$, are displayed in the top and bottom panels of Figure \ref{Fig15}, respectively. We can see that EKI not only estimates the shape (via the level-set function) of the regions with different conductivity but also the conductivity values on each region.
In Figure \ref{Fig15B} (top) we plot, as a function of $n$, the values of $\overline{L}_{1,f}$, $\overline{L}_{2,f}$, $\overline{\lambda}_{l}$, $\overline{\lambda}_{b}$ and $\overline{\lambda}_{h}$, i.e. the means of the ensembles $\{L_{1,f}^{(j)}\}_{j=1}^{J}$, $\{L_{2,f}^{(j)}\}_{j=1}^{J}$, $\{\lambda_{l}^{(j)}\}_{j=1}^{J}$, $\{\lambda_{b}^{(j)}\}_{j=1}^{J}$, and $\{\lambda_{h}^{(j)}\}_{j=1}^{J}$ (we reiterate that these variables are updated at each iteration of EKI). The ensemble mean for the intrinsic lengthscale of the leve-set function in the vertical direction is larger than the horizontal one. Again, this is due to the presence of the regions of low conductivity which are longer in the vertical direction.  From the middle-right panel we can see that the true conductivity in the background region $\overline{\lambda}_{b}^{\dagger}$ is recovered quite accurately. In contrast, the mean values $\overline{\lambda}_{l}$ and $\overline{\lambda}_{h}$ do not seem to vary much with respect to the mean. 

Although we are mainly focus on the deterministic case here, to further appreciate the accuracy of the inversion for these variables, in Figure \ref{Fig15B} (bottom) we show their probability densities approximated from the initial and final (converged) ensembles. For most of these variables we can see that the converged ensemble has a much smaller variance compared to the initial one. For the conductivity values we note that the lower and background values are identified accurately with the ensemble mean; the higher value is captured in the tail of the final ensemble. 


 \subsubsection{Comparisons of  \textbf{EKI-DMC} and  \textbf{EKI-LM}}

\begin{figure}[h!]
\centering
\includegraphics[scale=0.42]{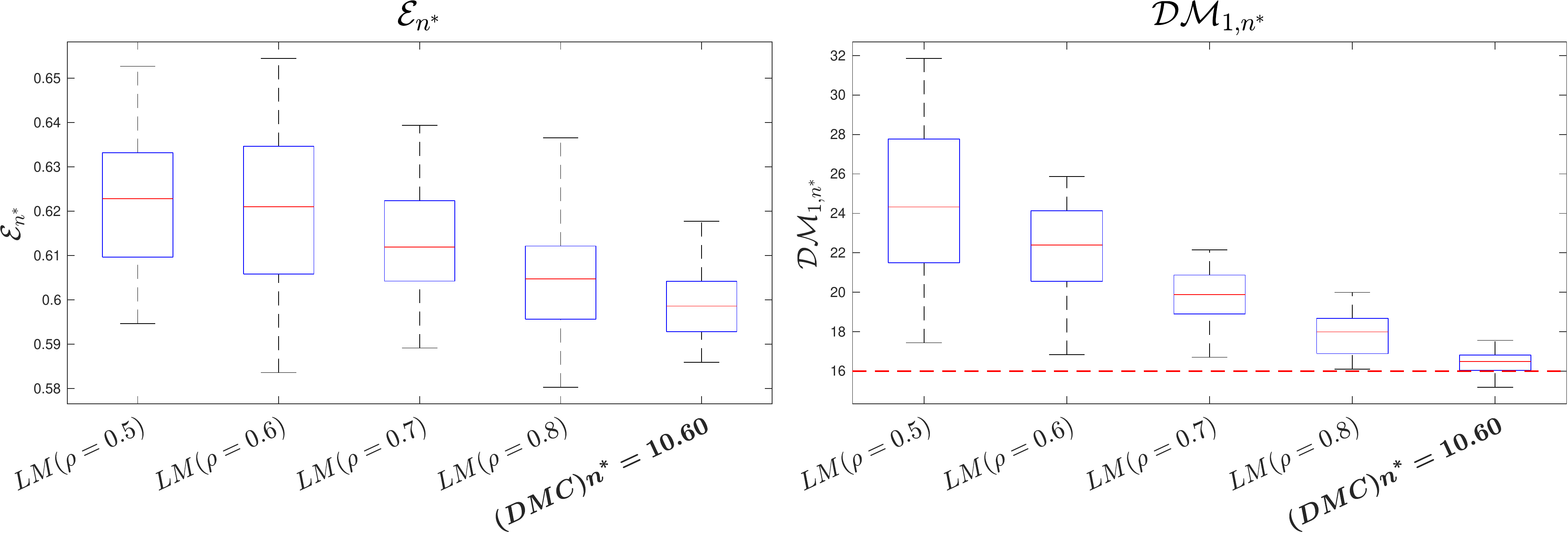}

\caption{{\bf Exp\_EIT$_2$.} Error with respect to the truth (left), $\mathcal{E}_{n^{*}}$ (see (\ref{eq:2003})), and data misfit $\mathcal{DM}_{1,n^*}$ (\ref{eq:2004}) (right) computed at the final iteration $n^{*}$ using {\bf EKI-DMC} and {\bf EKI-LM} with various choices of $\rho$. The noise level estimated by $\delta=\sqrt{M}$ is indicated with the dotted red-line in the right panel.}

\label{Fig12}
\end{figure}

\begin{figure}[h!]
\centering
\includegraphics[scale=0.45, trim=60 10 0 0, clip]{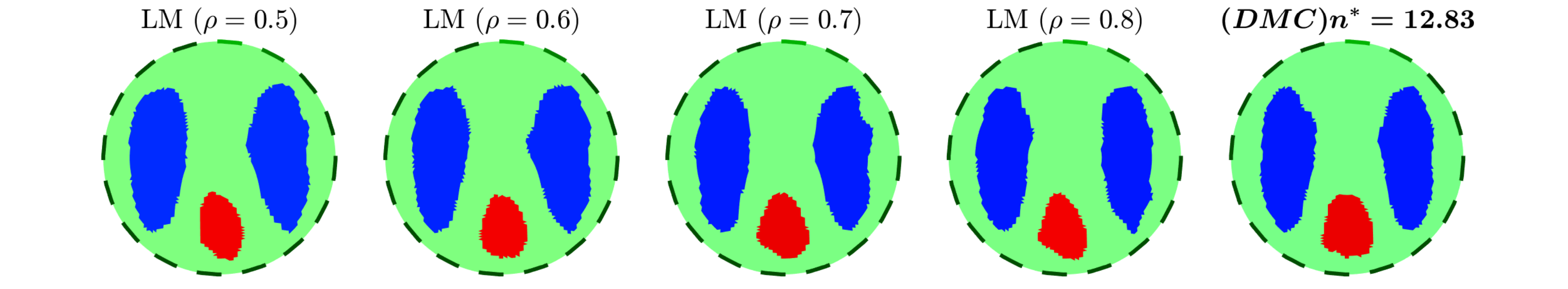} 
\caption{${\bf Exp\_EIT_2}$. Logarithm of $\kappa_{n^*}$ computed using the same initial ensemble ($J=200$) with {\bf EKI-LM} using different selections of the parameter $\rho$. In the right panel with display the log of $\kappa_{n^*}$ that we obtain using {\bf EKI-DMC}.}
\label{Fig13}
\end{figure}

In Figure \ref{Fig12} we compare the performance of \textbf{EKI-DMC} and  \textbf{EKI-LM} using 30 different initial ensembles. For \textbf{EKI-LM} we explore different choices of $\rho$. We note that \textbf{EKI-DMC} outperforms \textbf{EKI-LM} for all our choices of $\rho$. From Table \ref{Table1B} note that for $\rho>0.7$ the computational cost of  \textbf{EKI-LM} is approximately twice the cost of \textbf{EKI-DMC} and the cost even triple if we choose $\rho=0.8$. The plots of $\log{\kappa_{n}}$ from one run with the three algorithms (same initial ensemble) are shown in Figure \ref{Fig13}, where we see that all these runs perform well. Similarly conclusions to those in {\bf Exp\_EIT$_1$} are also drawn for this case. Namely, \textbf{EKI-DMC} is more accurate than \textbf{EKI-LM} in the chosen metrics although visually we achieve good performance from both algorithms and inputs. Experiments (not shown) indicate that substantial improvement in the performance of \textbf{EKI-LM} can be achieved by using the same stopping rule given by (\ref{eqA:9}) instead of (\ref{EQ2}). However, some choices of $\rho$ can result in inaccurate estimates in terms of error w.r.t the truth. The main advantage of \textbf{EKI-DMC} over \textbf{EKI-LM} is that the former does not rely on tuning parameters whose optimal selection is crucial for the stability and accuracy of the scheme.


\section{Conclusions}\label{Conclu}

We introduced the data misfit controller (DMC): a new adaptive regularisation strategy within the classical EKI setting. This led to an algorithm, {\bf EKI-DMC}, that in contrast to existing EKI approaches, it does not require any tuning parameters. Although we focus on the solution of deterministic identification problem, the proposed DMC is motivated from the Bayesian perspective of EKI within the tempering setting, where the inverse of the regularisation parameter, $\alpha_{n}^{-1}$, controls the transition between two consecutive intermediate measures. The Bayesian tempering setting provides us with a condition that these parameters must be satisfy ($\sum_{n=1}^{n^*}\alpha_{n}^{-1}=1$) to bridge the prior and posterior. We encode this condition for the termination of {\bf EKI-DMC} together with our new method for choosing $\alpha_{n}$. We show that the selection satisfies a heuristic statistical discrepancy principle which also controls Jeffreys' divergence between consecutive measures. 

We applied {\bf EKI-DMC} for the solution of EIT with the complete electrode model in which the conductivity was parameterised via different maps that enabled us, via the EKI framework applied to the appropriate parameters, to characterise both smooth and piece-wise constants conductivities. Whittle-Matern fields were at the core of these parameterisations which included the intrinsic lengthscales as inputs that we estimate within EKI. For the piece-wise constant case, we use a truncated WF field (the level-set function) to characterise discontinuities between different regions. As with any other EKI algorithms, our results show that the performance of {\bf EKI-DMC} relies on reasonable choices of ensemble size, $J$. For sufficiently large choices of $J$, our experiments show that {\bf EKI-DMC} is quite robust and capable of producing accurate identification of physical properties suitably parameterised.

We conducted a  performance comparison between {\bf EKI-DMC} and the {\bf EKI-LM} approach of \cite{EnsembleYo}. In most cases {\bf EKI-DMC} outperforms the accuracy of {\bf EKI-LM} in terms of error w.r.t the truth and data misfit, but we noted that suitable choices of tuning parameters and possible modifications to the stopping rule can produce comparable performance to {\bf EKI-DMC}. We recognise that similar performance can also be achieved by using different choices of $\alpha_{n}$ in EKI including those discussed in Section \ref{Liter}. We reiterate that, as with {\bf EKI-LM} and its variants, most regularisation approaches for EKI rely on the choice of additional tuning parameters which, of course, can be tuned to display comparable or even better performance to {\bf EKI-DMC}. However, there is no principled approach for the optimal selection of those parameters. Optimal tuning parameters can only be informed via careful numerical investigations on the given problem-specific setting. Our results show that {\bf EKI-DMC} is a robust self-tuning regularisation strategy of EKI, ideally suited for practical and operational settings for which finding optimal choices of tuning parameters in existing EKI approaches may not be computationally feasible. 


\section*{Acknowledgments}
This work was supported by the Engineering and Physical Sciences Research Council [grant number EP/P006701/1]; through the EPSRC Future Composites Manufacturing Research Hub.

  \bibliographystyle{elsarticle-num} 
  \bibliography{EKI_bib}

\appendix

\section{Motivation of EKI from the Bayesian tempering approach}\label{Ape_LM}
We introduce a series of approximations and Gaussian assumptions that lead to the EKI algorithm (Algorithm \ref{Al1}).

\subsection{Linearisation and Gaussian approximations}\label{Gauss_approx}

Suppose that the collection of tempering parameters $\{\phi_{n}\}_{n=1}^{N}$ in (\ref{eqA:6B}) have been specified. Our objective now is to construct a sequence of Gaussian approximations of each measure $\mu_{n}$ in (\ref{eqA:5}). To this end, let $\nu_{0}=N(m_0,\C_0)$ be a Gaussian approximation of the prior $\mu_{0}$, and let us denote by $D\cG$ the Frechet derivative of $\cG$.  We recursively construct a sequence of Gaussian approximations $\{\nu_{n}=N(m_n,\C_n)\}_{n=1}^{N+1}$ of $\{\mu_{n}\}_{n=1}^{N+1}$ via the  following expression
\begin{eqnarray}\label{eqA:11}
\nu_{n+1}(du)\propto \nu_{n}(du) \exp\Big[-\frac{1}{2}\norm{ (\alpha_{n}\Gamma)^{-1/2}(y-\cG_n-D\cG_{n}(u-m_{n}))}^2\Big],
\end{eqnarray}
where for ease in the notation we have defined $\cG_{n}\equiv \cG(m_{n})$ and $D\cG_{n}\equiv D\cG(m_{n})$. We note the right hand side of (\ref{eqA:11}) involves the linearisation of the forward map around the mean of $\nu_{n}=N(m_n,\C_n)$.

Recursive formulas for the mean and covariance of $\nu_n=N(m_n,\C_n)$ can be obtained by completing the square in (\ref{eqA:11}) (see Theorem 6.20 in \cite{Andrew}). Indeed, since $\nu_{n}=N(m_n,\C_n)$ and the model $L(u)\equiv D\cG_{n}(u-m_{n})$ is linear,  then $\nu_{n+1}=N(m_{n+1},\C_{n+1})$ with 
\begin{eqnarray}\label{eqA:12}
m_{n+1}&=&m_{n}+\C_{n}D\cG_n^{*}(D\cG_n\C_{n}\,D\cG_n^{*} +\alpha_{n}\Gamma   )^{-1}(y-\cG_n)\\
\C_{n+1}&=&\C_{n}-\C_{n}D\cG_n^{*}(D\cG_n\C_{n}\,D\cG_n^{*} +\alpha_{n}\Gamma   )^{-1}D\cG_n\C_{n}\label{eqA:13}
\end{eqnarray}
where $D\cG_n^{*}$ denotes the adjoint of $D\cG$ evaluated at $m_{n}$. 

\begin{remark}[Levenberg-Marquardt from linearised Bayesian tempering]
Using standard arguments (see for example Lemma 3.1 in \cite{EnsembleYo}) it can be shown that $m_{n}$ in (\ref{eqA:12}) satisfies (\ref{eqA:14}). In the case where $\C_{n}$ is the identity operator, (\ref{eqA:14}) yields the standard Levenberg-Marquardt (LM) iterative scheme \cite{Mar}. For the modified version in (\ref{eqA:14}), we note that the recursive formula for the mean involves introducing the precision operator for $u$ in the regularisation term in the right hand side of (\ref{eqA:14}). 
\end{remark}
\begin{remark}[The linear-Gaussian case]
Note that, if $\mu_{0}=\nu_{0}$ (i.e. the prior is Gaussian) and the forward map $\cG$ is linear, then we have that $D\cG_{n}(u-m_{n})=\cG(u-m_{n})=\cG(u)-\cG_{n}$. Hence, (\ref{eqA:11}) and (\ref{eqA:7}) coincide and so $\mu_{n}=\nu_{n}$ for all $n=0,\dots, N+1$. In particular, in the linear-Gaussian case the final measure $\nu_{n}^{*}$ coincides with the posterior (i.e. $\nu_{N+1}=\mu_{N+1}=\mu)$
\end{remark}

\subsection{Derivative-free ensemble approximation.}\label{derivative_free}

We now introduce further approximations that will lead to the EKI algorithm under consideration. Let us denote by $u_{n}$ a random variable such that $u_{n}\sim \nu_{n}=N(m_{n},\C_{n})$. Denote by $\E_{n}$ expectation with respect to $\nu_{n}$. Let us consider the first order approximation
\begin{eqnarray}\label{eqA:14B}
\cG(u_{n})\simeq \cG_{n}+D\cG_{n}(u_{n}-m_{n})
\end{eqnarray}
that we used in (\ref{eqA:11}) to define our Gaussian approximations for the tempering scheme introduced in \ref{Gauss_approx}. From (\ref{eqA:14B}) it follows that 
\begin{eqnarray}\label{eqA:15}
\E_{n}[\cG(u_{n})]\simeq \cG_n,\qquad \cG(u_n)-\E_{n}[\cG(u_{n})]\simeq D\cG_{n}(u_n-m_{n})
\end{eqnarray}
Hence,
\begin{eqnarray}\label{eqA:16}
  \co_{n}(u_n,\cG(u_{n}))\equiv &\E_{n}[( u_{n}-m_{n})\otimes (\cG(u_n)-\E_{n}[\cG(u_{n})]) ]&\simeq  \C_{n}D\cG_n^{*},\\
\co_{n}(\cG(u_{n}))\equiv &\E_{n}[ (\cG(u_n)-\E_{n}[\cG(u_{n})] )\otimes (\cG(u_n)-\E_{n}[\cG(u_{n})] )]&\simeq D\cG_n\C_{n}D\cG_n^{*}\label{eqA:16B}
\end{eqnarray}
where we have used the fact that $\E_{n}[( u_{n}-m_{n})\otimes ( u_{n}-m_{n})] =\C_{n}$. If we use approximations (\ref{eqA:15})-(\ref{eqA:16B}) in (\ref{eqA:12})-(\ref{eqA:13}) we note that
\begin{eqnarray}\label{eqA:17}
m_{n+1}&\simeq &\tilde{m}_{n+1}\equiv m_{n}+\co_{n}(u_n,\cG(u_n))(\co_{n}(\cG(u_{n}))+\alpha_{n}\Gamma   )^{-1}(y-\E_{n}[\cG(u_{n})])\\
\C_{n+1}&\simeq&\tilde{\C}_{n+1}\equiv \C_{n}-\co_{n}(u_n,\cG(u_n))(\co_{n}(\cG(u_{n}))+\alpha_{n}\Gamma   )^{-1}\co_{n}(\cG(u_n),u_n)\label{eqA:18}
\end{eqnarray}
These approximations to the mean and covariance of the sequence of approximate measures $\{\nu_{n}\}_{n=1}^{N+1}$ do not involve derivatives of the forward map. However, the covariance and cross-covariance that appear in (\ref{eqA:17})-(\ref{eqA:18}) cannot be computed in closed form. This issue is overcome by using particle approximations of the each approximate measure $\nu_{n}=N(m_{n},\C_{n})$. In other words, we consider approximations
\begin{eqnarray}\label{eqnewA1}
\nu_{n}^{J}(u_{n})=\frac{1}{J}\sum_{j=1}^{J}\delta(u_{n}-u_{n}^{(j)})
\end{eqnarray}
where we assume that $u_{n}^{(j)}\sim N(m_{n},\mathcal{C}_{n})$. The classical EKI update formula (\ref{equ2}) for the ensemble of particles $\{u_{n}^{(j)}\}_{j=1}^{J}$ is defined in such a way, that the corresponding ensemble mean and covariance  approximate those in (\ref{eqA:17})-(\ref{eqA:18}) as $J\to \infty$. To see this more clearly, let us note that the ensemble approximation of $\E_{n}[u_{n}]$ and $\E_{n}[\cG(u_{n})]$ are $\overline{u}_{n}$ and $\overline{\cG}_{n}$ defined in (\ref{equ4}). The ensemble approximations of $\co_{n}(\cG(u_{n}))$, $\co_{n}(u_n,\cG(u_{n}))$ and $\C_{n}$, denoted by $C_{n}^{\cG\cG}$, $C_{n}^{u\cG}$ and $\C_{n}^{uu}$, are defined by (\ref{equ3}), (\ref{equ3B}) and 
\begin{eqnarray}\label{eqnewA2}
\C_{n}^{uu} \equiv \frac{1}{J-1}\sum_{j=1}^{J} (u_{n}^{(j)}-\overline{u}_{n})\otimes(u_{n}^{(j)}-\overline{u}_{n}),
\end{eqnarray}
respectively. From (\ref{equ2}) we note that
\begin{eqnarray}\label{eq:m16}
\ou_{n+1}&=&\ou_{n}+\C_{n}^{u\cG}(\C_{n}^{\cG \cG} +\alpha_{n}\Gamma  )^{-1}(y+\sqrt{\alpha_n}\overline{\xi}_{n}-\overline{\cG}_{n})
\end{eqnarray}
It can be shown  (see for example \cite{Mandel} for a rigorous proof in finite dimensions) that 
\begin{eqnarray}\label{eq:m16B}
\ou_{n+1}\to \tilde{m}_{n+1}\equiv m_{n}+\co_{n}(u_n,\cG(u_n))(\co_{n}(\cG(u_{n}))+\alpha_{n}\Gamma   )^{-1}(y-\E_{n}[\cG(u_{n})])\simeq m_{n+1}\\
\C_{n+1}^{uu}\to \tilde{\C}_{n+1}\equiv \C_{n}-\co_{n}(u_n,\cG(u_n))(\co_{n}(\cG(u_{n}))+\alpha_{n}\Gamma   )^{-1}\co_{n}(\cG(u_n),u_n)\simeq \C_{n+1}.\label{eq:m16b}
\end{eqnarray}
as $J\to \infty$. Moreover, for the particles in (\ref{equ2}) we have
$$\nu_{n+1}^{J}(u_{n+1})=\frac{1}{J}\sum_{j=1}^{J}\delta(u_{n+1}-u_{n+1}^{(j)})\to \tilde{\nu}_{n+1}\equiv N(\tilde{m}_{n+1},\tilde{\C}_{n+1}).$$

The development above (informally) shows that if $\nu_{n}^{J}\simeq \nu_{n}=N(m_{n},\C_{n})$, then $\tilde{\nu}_{n+1}^{J}\simeq \tilde{\nu}_{n+1}= N(\tilde{m}_{n+1},\tilde{\C}_{n+1})$. Furthermore, if (\ref{eqA:15})-(\ref{eqA:16B}) are accurate enough, then $N(\tilde{m}_{n+1},\tilde{\C}_{n+1})\simeq \nu_{n+1}=N(m_{n+1},\C_{n+1})$ and so the EKI ensemble $\tilde{\nu}_{n+1}^{J}$ approximates $\nu_{n+1}$. Recall that the $\nu_{n}$'s are Gaussian approximations of the tempered distributions $\mu_{n}$. Hence the regularisation parameter $\alpha_{n}$ in EKI is the inverse of the difference between consecutive tempering parameters (see eq. (\ref{eqA:8})) which is, in turn, a Gaussian approximation of $\mu_{n+1}$.


\begin{remark}[Squared-root EKI]
It is worth noticing that other approaches can be used to approximate $N(\tilde{m}_{n+1},\tilde{\C}_{n+1})$ above. This includes the so-called ensemble square-root formulations \cite{doi:10.1137/140965363} in which the particles are cleverly updated so that their sample mean and covariance coincide (exactly) with $\tilde{m}_{n+1}$ and $\tilde{\C}_{n+1}$. While these has been shown to be beneficial for very small samples (i.e. $<50$), we note that $N(\tilde{m}_{n+1},\tilde{\C}_{n+1})$ does not, in general, coincides with $\nu_{n+1}=N(m_{n+1},\C_{n+1})$ (unless $\cG$ is linear).
\end{remark}

\begin{remark}[EKI as a derivative-free approximation of LM]\label{LM}
For sufficiently large $J$, the mean of the ensemble $\ou_{n+1}$ approximates $\tilde{m}_{n+1}$ and so $m_{n+1}$ which is, in turn, the iteration of the LM scheme in (\ref{eqA:14}) (see Remark \ref{LM}). Therefore, we can interpret EKI as a derivative-free approximation of the LM scheme constrained to the subspace generated by the initial ensemble $\{u_{0}^{(j)}\}_{j=1}^{J}$. More specifically, the ensemble mean $\overline{u}_{n}$ define by the recursive formula (\ref{equ4}) satisfies the following subspace invariance property (see Theorem 2.1 in \cite{ILS13a})
$$\overline{u}_{n}\in \mathcal{S}_{0}\equiv \text{span}\{u_{0}^{(j)}\}_{j=1}^{J}$$
for all $n\in \mathbb{N}$. We expect EKI to produce approximate solutions to (\ref{equ2}) within the subspace defined above. While the numerical experiments from \cite{EnsembleYo} provides evidence of such a claim, to the best of our knowledge, the convergence of EKI in this context is still an open problem.
\end{remark}

\section{Lemmas, Theorems, and Poofs} \label{app: appendix}

Let $y\in \mathbb{R}^{M}$ fixed. We recall that from Assumption \ref{ass: forward map} on the forward map $\cG$, the data misfit $\Phi(\cdot ,y): \mathcal{H} \to [0,+\infty)$ is a square integrable functional on the probability space $(\mathcal{H}, \mathcal{B}(\mathcal{H}), \mu_0 )$, i.e. 
\begin{equation} \label{eq: condition TPM NZ 2}
\int_\mathcal{H} |\Phi(u;y)|^2 \,\mu_0(\mathrm{d}u) < \infty
\end{equation}
This implies that $\Phi(\cdot ,y): \mathcal{H} \to [0,+\infty)$ is absolutely integrable functional on the probability space $(\mathcal{H}, \mathcal{B}(\mathcal{H}) , \mu_0 )$, i.e. 
\begin{equation} \label{eq: condition TPM NZ 1}
\int_\mathcal{H} |\Phi(u;y)| \,\mu_0(\mathrm{d}u) \leq M_0 < \infty
\end{equation}

\subsection{Proof of Lemma \ref{lem: equivalence of TPM}}\label{proof_lema}
\begin{proof}
We use definition \eqref{eq: define N_t},  Jensen's inequality and \eqref{eq: condition TPM NZ 1} to see that 
\begin{align*}
N_t & = \int_\mathcal{H} \exp(-t\Phi(u;y)) \, \mu_0(\mathrm{d}u)  \geq \exp\left( -  t\int_\mathcal{H} \Phi(u;y) \, \mu_0(\mathrm{d}u) \right) \geq \exp\left( -  tM_0 \right) > 0,
\end{align*}
which shows that $N_{t}$ is strictly positive for all bounded $t\geq 0$, so the right hand side of formula \eqref{eq: Radon-Nikodym derivative mu_t mu_0} is well-defined. By Definition \eqref{eq: Radon-Nikodym derivative mu_t mu_0}, it then follows that $\mu_t \ll \mu_0$. We now prove $\mu_0 \ll \mu_t$, i.e. $\forall \mathcal{X} \in \mathcal{B}(\mathcal{H})$, $u_t(\mathcal{X})=0 \implies u_0(\mathcal{X})=0 $. We prove this statement by contradiction. Assume that there exists $\mathcal{X} \in \mathcal{B}(\mathcal{H})$ such that $u_t(\mathcal{X})=0$ and $u_0(\mathcal{X}) > 0$. Since $u_0(\mathcal{X}) > 0$, we can define a another probability measure, $\nu_\mathcal{\mathcal{X}}$, on the measurable space $(\mathcal{X},\mathcal{B}(\mathcal{X}))$, such that for all $\mathcal{S}\in \mathcal{B}(\mathcal{X})$,
\begin{equation} \label{eq: X measure change}
\nu_\mathcal{\mathcal{X}}(\mathcal{S}):=\frac{u_0(\mathcal{S})}{u_0(\mathcal{X})}
\end{equation}
Using $\mu_t \ll \mu_0$ and expressions \eqref{eq: Radon-Nikodym derivative mu_t mu_0} and \eqref{eq: X measure change} we have 
\begin{align*}
u_t(\mathcal{X}):
&=     \int_\mathcal{X} \,\mu_t(du) 
=     \int_\mathcal{X} \frac{\mathrm{d}\mu_t}{\mathrm{d}\mu_0}(u) \,\mu_0(du) 
=     \frac{1}{N_t} \int_\mathcal{X} \exp(-t\Phi(u;y)) \,\mu_0(du) \\
=    & \frac{\mu_0(\mathcal{X})}{N_t} \left( \int_\mathcal{X} \exp(-t\Phi(u;y))\, \nu_\mathcal{\mathcal{X}}(\mathrm{d}u) \right) 
\end{align*}
We use Jensen's inequality, formula \eqref{eq: X measure change} as well as \eqref{eq: condition TPM NZ 1} to obtain
\begin{align*}
u_t(\mathcal{X})\geq & \frac{\mu_0(\mathcal{X})}{N_t}   \exp\left( -t\int_\mathcal{X} \Phi(u;y) \, \nu_\mathcal{\mathcal{X}}(\mathrm{d}u) \right)
=     \frac{\mu_0(\mathcal{X})}{N_t}   \exp\left( - \frac{t}{\mu_0(\mathcal{X})} \int_\mathcal{X} \Phi(u;y) \, \mu_0(\mathrm{d}u) \right) \\
\geq & \frac{\mu_0(\mathcal{X})}{N_t}   \exp\left( -  \frac{t  M_0}{\mu_0(\mathcal{X})} \right) > 0  
\end{align*}
which contradicts our assumption that $u_t(\mathcal{X})=0$. 
\end{proof}

\subsection{Proof of Theorem \ref{the: dynamic equation}}\label{proof_of_theo}
Let us first prove:
\begin{theorem} \label{the: derivative NZ}
For any bounded $t \geq 0$, the normalizing constant $N_t$ defined in formula \eqref{eq: define N_t}, as a function of $t$, is differentiable and its derivative $N_t'$ is expressed by
\begin{equation} \label{eq: derivative NZ}
\frac{N_t'}{N_t} = - \int_\mathcal{H} \Phi(u;y) \mu_t(\mathrm{d}u)
\end{equation}
where $\mu_t$ is the probability measure determined via formula \eqref{eq: Radon-Nikodym derivative mu_t mu_0}.
\begin{proof}
For any $u \in \mathcal{H}$ and for any bounded $t\geq 0$, let us define 
\begin{equation*}
f_u(t) := \exp(-t\Phi(u;y))
\end{equation*}
Then, $N_t$ defined in formula \eqref{eq: define N_t} can be rewritten as
\begin{equation*}
N_t = \int_\mathcal{H} f_u(t) \, \mu_0(\mathrm{d}u)
\end{equation*}
Notice that the derivative of $f_u$ is bounded by $\Phi(u;y)$ for any $t \geq 0$, since
\begin{equation*}
|f_u'(t)| = |-\Phi(u)\exp(-t\Phi(u;y))| \leq \Phi(u)
\end{equation*}
Recall that $\Phi(\cdot;y)$ is absolutely integrable. Thus, according to the dominated convergence theorem, the derivative of $N_t$ can be calculated by
\begin{equation*}
N_t' = \left( \int_\mathcal{H} f_u(t) \, \mu_0(\mathrm{d}u) \right)' = \int_\mathcal{H} f_u'(t) \, \mu_0(\mathrm{d}u) = - \int_\mathcal{H} \Phi(u)\exp(-t\Phi(u)) \, \mu_0(\mathrm{d}u)
\end{equation*}
From the proof of Lemma \ref{lem: equivalence of TPM} we know $N_t$ is strictly greater than $0$ for any bounded $t \geq 0$. Therefore, $N_t'$ can be divided by $N_t$,
\begin{equation*}
\frac{N_t'}{N_t} = - \frac{1}{N_t} \int_\mathcal{H} \Phi(u)\exp(-t\Phi(u;y)) \, \mu_0(\mathrm{d}u)
\end{equation*}
Expression \eqref{eq: derivative NZ} follows from Lemma \ref{lem: equivalence of TPM} which establishes the equivalence between $\mu_t$ and $\mu_0$.
\end{proof}
\end{theorem}

We now prove Theorem \ref{the: dynamic equation}:
\begin{proof}
$\mathbb{E}\left\lbrace g(u_t) \right\rbrace $ can be rewritten by changing the measure from $\mu_t$ to $\mu_0$ (recall $\mu_t \ll \mu_0$),
\begin{align*}
\mathbb{E}\left\lbrace g(u_t) \right\rbrace := & \int_\mathcal{H} g(u)\, \mu_t(\mathrm{d}u) = \int_\mathcal{H} g(u)\frac{\mathrm{d}\mu_t}{\mathrm{d}\mu_0}(u) \, \mu_0(\mathrm{d}u) \\
= & \frac{1}{N_t} \int_\mathcal{H} g(u)\exp(-t\Phi(u)) \, \mu_0(\mathrm{d}u) 
=  \frac{1}{N_t} \int_\mathcal{H} f_u(t) \, \mu_0(\mathrm{d}u)
\end{align*}
where $f_u(t)$ is defined by 
\begin{equation*}
f_u(t) := g(u) \exp(-t\Phi(u))
\end{equation*}
for any $u \in \mathcal{H}$ and any bounded $t \geq 0$. By applying the chain rule, the derivative of $\mathbb{E}\left\lbrace g(u_t) \right\rbrace$ is given by
\begin{align}
\mathbb{E}\left\lbrace g(u_t) \right\rbrace' & = \left(  \frac{1}{N_t} \int_\mathcal{H} f_u(t) \, \mu_0(\mathrm{d}u)  \right)' 
 = -  \frac{ N_t'\int_\mathcal{H} f_u(t) \, \mu_0(\mathrm{d}u)}{N_t^2} + \frac{\left( \int_\mathcal{H} f_u(t) \, \mu_0(\mathrm{d}u) \right)'}{N_t} \label{eq: dynamic equation__ 1}
\end{align}
Notice that the derivative of $f_u$ is bounded by $|g(u)\Phi(u)|$ for any $t \geq 0$, since
\begin{equation*}
\left| f_u'(t) \right| =  |-g(u)\Phi(u)\exp(-t\Phi(u))| \leq |g(u)\Phi(u)|
\end{equation*}
Since both $g$ and $\Phi$ are square integrable, it follows from Cauchy-Schwarz inequality that $g\cdot\Phi$ is absolutely integrable. Thus, according to the dominated convergence theorem, we have
\begin{equation} \label{eq: dynamic equation__ 2}
\left( \int_\mathcal{H} f_u(t) \, \mu_0(\mathrm{d}u) \right)' = \int_\mathcal{H} f_u'(t) \, \mu_0(\mathrm{d}u) = -\int_\mathcal{H} \Phi(u)f_u(t) \, \mu_0(\mathrm{d}u)
\end{equation}
On the other hand, we note that the conditions of Theorem \ref{the: derivative NZ} applies.  Then, we substitute \eqref{eq: dynamic equation__ 2} and \eqref{eq: derivative NZ} in \eqref{eq: dynamic equation__ 1} to arrive at
\begin{equation*}
\mathbb{E}\left\lbrace g(u_t) \right\rbrace' 
= \frac{\int_\mathcal{H} \Phi(u) \mu_t(\mathrm{d}u) \int_\mathcal{H} f_u(t) \, \mu_0(\mathrm{d}u)}{N_t} -\frac{\int_\mathcal{H} \Phi(u)f_u(t) \, \mu_0(\mathrm{d}u) }{N_t}
\end{equation*}
According to Lemma \ref{lem: equivalence of TPM}, $\mu_t$ and $\mu_0$ are equivalent for any bounded $t \geq 0$. Hence, we change measure in the right hand side of the equation above using formula \eqref{eq: Radon-Nikodym derivative mu_t mu_0}. Thus, 
\begin{equation*}
\mathbb{E}\left\lbrace g(u_t) \right\rbrace' 
= \int_\mathcal{H} \Phi(u) \mu_t(\mathrm{d}u) \int_\mathcal{H} g(u) \, \mu_t(\mathrm{d}u) - \int_\mathcal{H} \Phi(u)g(u) \, \mu_t(\mathrm{d}u)
\end{equation*}
Since both $\Phi$ and $g$ are square integrable, it is not difficult to see that the right hand side in the previous expression gives the desired result in \eqref{eq: dynamic equation}.
\end{proof}

\end{document}